\documentstyle[12pt,fleqn]{article}

\begin{document}

\title{Complex Numbers in Three Dimensions}

\author{Silviu Olariu
\thanks{e-mail: olariu@ifin.nipne.ro}\\
Institute of Physics and Nuclear Engineering, Tandem Laboratory\\
76900 Magurele, P.O. Box MG-6, Bucharest, Romania}

\date{4 August 2000}

\maketitle

\abstract

A system of commutative hypercomplex numbers of the form $w=x+hy+kz$ are
introduced in 3 dimensions, the variables $x, y$ and $z$ being real numbers.
The multiplication rules for the complex units $h, k$ are $h^2=k, k^2=h, hk=1$.
The operations of addition and multiplication of the tricomplex numbers
introduced in this paper have a simple geometric interpretation based on the
modulus $d$, amplitude $\rho$, polar angle $\theta$ and azimuthal angle $\phi$.
Exponential and trigonometric forms are obtained for the tricomplex numbers,
depending on the variables $d$, $\rho$, $\theta$ and $\phi$. The tricomplex
functions defined by series of powers are analytic, and the partial derivatives
of the components of the tricomplex functions are closely related. The
integrals of tricomplex functions are independent of path in regions where the
functions are regular. The fact that the exponential form of the tricomplex
numbers contains the cyclic variable $\phi$ leads to the concepts of pole and
residue for integrals of tricomplex functions on closed paths. The polynomials
of tricomplex variables can be written as products of linear or quadratic
factors.

\endabstract

\section{Introduction}

A regular, two-dimensional complex number $x+iy$ 
can be represented geometrically by the modulus $\rho=(x^2+y^2)^{1/2}$ and 
by the polar angle $\theta=\arctan(y/x)$. The modulus $\rho$ is multiplicative
and the polar angle $\theta$ is additive upon the multiplication of ordinary 
complex numbers.

The quaternions of Hamilton are a system of hypercomplex numbers
defined in four dimensions, the
multiplication being a noncommutative operation, \cite{1} 
and many other hypercomplex systems are
possible, \cite{2a}-\cite{2b} but these hypercomplex systems 
do not have all the required properties of regular, 
two-dimensional complex numbers which rendered possible the development of the 
theory of functions of a complex variable.

A system of hypercomplex numbers in three dimensions is described in this work,
for which the multiplication is associative and commutative, and which is 
rich enough in properties so that exponential and trigonometric forms exist
for these numbers, and the concepts
of analytic tricomplex 
function,  contour integration and residue can be defined.
The tricomplex numbers introduced in this work have 
the form $u=x+hy+kz$, the variables $x, y$ and $z$ being real 
numbers.  
The multiplication rules for the complex units $h, k$ are
$h^2=k, \:  k^2=h,\:  hk=1$. In a geometric representation, the
tricomplex number $u$ is represented by the point $P$ of coordinates $(x,y,z)$.
If $O$ is the origin of the $x,y,z$
axes, $(t)$ the trisector line $x=y=z$ of the positive octant and 
$\Pi$ the plane $x+y+z=0$ passing through the origin $(O)$ and
perpendicular to $(t)$, then the tricomplex number $u$ can be described by
the projection $s$ of the segment $OP$ along the line $(t)$, by the distance
$D$ from $P$ to the line $(t)$, and by the azimuthal angle $\phi$ of the
projection of $P$ on the plane $\Pi$, measured from an angular origin defined
by the intersection of the plane determined by the line $(t)$ and the x axis,
with the plane $\Pi$. 
The amplitude $\rho$ of a twocomplex number is defined as 
$\rho=(x^3+y^3+z^3-3xyz)^{1/3}$, the polar 
angle $\theta$ of $OP$ with respect to the trisector line $(t)$ is given by
$\tan\theta=D/s$, and $d^2=x^2+y^2+z^2$.  
The amplitude $\rho$ is equal to
zero on the trisector line $(t)$ and on the plane $\Pi$. The division 
$1/(x+hy+kz) $is possible provided that $\rho\not= 0$.
The product of two tricomplex numbers is equal to zero if both numbers are
equal to zero, or if one of
the tricomplex numbers lies in the $\Pi$ plane and the other on the $(t)$
line. 

If $u_1=x_1+hy_1+kz_1, 
u_2=x_2+hy_2+kz_2$ are tricomplex
numbers of amplitudes and angles $\rho_1,\theta_1,\phi_1$ and respectively
$\rho_2, \theta_2, \phi_2$, then the amplitude and the angles $\rho,
\theta, \phi$ for
the product tricomplex number
$u_1u_2=x_1x_2+y_1z_2+y_2z_1+h(z_1z_2+x_1y_2+y_1x_2)+k(y_1y_2+x_1z_2+z_1x_2)$
are $\rho=\rho_1\rho_2,
\tan\theta=\tan\theta_1\tan\theta_2/\sqrt{2}, 
\phi=\phi_1+\phi_2$. 
Thus, the amplitude $\rho$ and  $(\tan\theta)/\sqrt{2}$ are
multiplicative quantities  
and the angle $\phi$ is an additive quantity upon the
multiplication of tricomplex numbers, which reminds the properties of 
ordinary, two-dimensional complex numbers.

For the description of the exponential function of a tricomplex variable, it
is useful to define the cosexponential functions 
${\rm cx}(\xi)=1+\xi^3/3!+\xi^6/6!\cdots, 
{\rm mx}(\xi)=\xi+\xi^4/4!+\xi^7/7!\cdots, 
{\rm px}(\xi)=\xi^2/2+\xi^5/5!+\xi^8/8!\cdots $, where p and m stand for plus
and respectively minus, as a reference to the sign of a phase shift in the
expressions of these functions. These functions fulfil the relation ${\rm cx}^3
\xi +{\rm px}^3 \xi+{\rm mx}^3 \xi -3 {\rm cx} \xi \:
{\rm px} \xi\: {\rm mx} \xi =1$. 

The exponential form of a tricomplex number is
$u=\rho$
$\exp\left[(1/3)(h+k)\ln(\sqrt{2}/\tan\theta)\right.$
$\left.+(1/3)(h-k)\phi\right]$,
and the trigonometric form of the tricomplex number is
$u=d\sqrt{3/2}$
$\left\{(1/3)(2-h-k)\sin\theta+\right.$
$\left.(1/3)(1+h+k)\sqrt{2}\cos\theta\right\}$
$\exp\left\{(h-k)\phi/\sqrt{3}\right\}$. 

Expressions are given for the elementary functions of tricomplex variable.
Moreover, it is shown that the region of convergence of series of powers of
tricomplex variables are cylinders with the axis parallel to the trisector
line. 
A function $f(u)$ of the tricomplex variable 
$u=x+hy+kz$ can be defined by a corresponding 
power series. It will be shown that the function $f(u)$ has a
derivative at $u_0$ independent of the
direction of approach of $u$ to $u_0$. If the tricomplex function $f(u)$
of the tricomplex variable $u$ is written in terms of 
the real functions $F(x,y,z),G(x,y,z),H(x,y,z)$ of real variables $x,y,z$ as
$f(u)=F(x,y,z)+hG(x,y,z)+kH(x,y,z)$, then relations of equality 
exist between partial derivatives of the functions $F,G,H$, and the differences
$F-G, F-H, G-H$ are solutions of the equation of Laplace.

It will be shown that the integral $\int_A^B f(u) du$ of a regular tricomplex
function between two points $A,B$ is independent of the three-dimensional
path connecting the points $A,B$.
If $f(u)$ is an analytic tricomplex function, then $\oint_\Gamma f(u)
du/(u-u_0) = 
2\pi ( h-k) f(u_0)$ if the integration loop is threaded by the parallel through
$u_0$ to the line $(t)$.

A tricomplex polynomial  $u^m+a_1 u^{m-1}+\cdots+a_{m-1} u +a_m $ can be
written as a 
product of linear or quadratic factors, although the factorization may not be
unique. 

This paper belongs to a series of studies on commutative complex numbers in $n$
dimensions. \cite{2c}
The tricomplex numbers described in this work are a particular case for 
$n=3$ of the polar hypercomplex numbers in $n$ dimensions.\cite{2c},\cite{2d}

\section{Operations with tricomplex numbers}

A tricomplex number is determined by its three components $(x,y,z)$. The sum
of the tricomplex numbers $(x,y,z)$ and $(x^\prime,y^\prime,z^\prime)$ is the
tricomplex 
number $(x+x^\prime,y+y^\prime,z+z^\prime)$. The product of the tricomplex
numbers 
$(x,y,z)$ and $(x^\prime,y^\prime,z^\prime)$ is defined in this work to be the
tricomplex 
number
$(xx^\prime+yz^\prime+zy^\prime,zz^\prime+xy^\prime+yx^\prime,
yy^\prime+xz^\prime+zx^\prime)$.

Tricomplex numbers and their operations can be represented by writing the
tricomplex number $(x,y,z)$ as  
$u=x+hy+kz$, where $h$ and $k$ are bases for which the
multiplication rules are 
\begin{equation}
 h^2=k, \:  k^2=h, \: 1\cdot h=h,\:  1\cdot k =k, \:  hk=1. 
\label{1}
\end{equation}
Two tricomplex numbers $u=x+hy+kz, u^\prime=x^\prime+hy^\prime+
kz^\prime$ are equal, $u=u^\prime$, if and only if $x=x^\prime, y=y^\prime,
z=z^\prime$. 
If $u=x+hy+kz, 
u^\prime=x^\prime+hy^\prime+kz^\prime$ are tricomplex
numbers, the sum $u+u^\prime$ and the 
product $uu^\prime$ defined above can be obtained by applying the usual
algebraic rules to the sum 
$(x+hy+kz)+(x^\prime+hy^\prime+kz^\prime)$
and to the product 
$(x+hy+kz)(x^\prime+hy^\prime+kz^\prime)$,
and grouping of the resulting terms,
\begin{equation}
u+u^\prime=x+x^\prime+h(y+y^\prime)+k(z+z^\prime),
\label{1a}
\end{equation}
\begin{equation}
uu^\prime=xx^\prime+yz^\prime+zy^\prime+h(zz^\prime
+xy^\prime+yx^\prime)+k(yy^\prime+xz^\prime+zx^\prime) .
\label{1b}
\end{equation}
If $u,u^\prime,u^{\prime\prime}$ are tricomplex numbers, the multiplication is
associative 
\begin{equation}
(uu^\prime)u^{\prime\prime}=u(u^\prime u^{\prime\prime})
\label{2}
\end{equation}
and commutative
\begin{equation}
u u^\prime=u^\prime u ,
\label{3}
\end{equation}
as can be checked through direct calculation.
The tricomplex zero is $0+h\cdot 0+k\cdot 0,$ denoted simply 0, 
and the tricomplex unity is $1+h\cdot 0+k\cdot 0,$ 
denoted simply 1.

The inverse of the tricomplex number $u=x+hy+kz$ is a
tricomplex number $u^\prime=x^\prime+y^\prime+z^\prime$ having the property
that 
\begin{equation}
uu^\prime=1 .
\label{4}
\end{equation}
Written on components, the condition, Eq. (\ref{4}), is
\begin{equation}
\begin{array}{c}
xx^\prime+zy^\prime+yz^\prime=1 ,\\
yx^\prime+xy^\prime+zz^\prime=0 ,\\
zx^\prime+yy^\prime+xz^\prime=0 .
\end{array}
\label{5}
\end{equation}
The system (\ref{5}) has the solution
\begin{equation}
x^\prime=\frac{x^2-yz}{x^3+y^3+z^3-3xyz} ,
\label{6a}
\end{equation}
\begin{equation}
y^\prime=\frac{z^2-xy}{x^3+y^3+z^3-3xyz} ,
\label{6b}
\end{equation}
\begin{equation}
z^\prime=\frac{y^2-xz}{x^3+y^3+z^3-3xyz} ,
\label{6c}
\end{equation}
provided that $x^3+y^3+z^3-3xyz\not=0$. Since
\begin{equation}
x^3+y^3+z^3-3xyz=(x+y+z)(x^2+y^2+z^2-xy-xz-yz) ,
\label{7}
\end{equation}
a tricomplex number $x+hy+kz$ has an inverse, unless
\begin{equation}
x+y+z=0  
\label{8}
\end{equation}
or
\begin{equation}
x^2+y^2+z^2-xy-xz-yz=0 . 
\label{9}
\end{equation}

The relation in Eq. (\ref{8}) represents the plane $\Pi$ perpendicular
to the trisector line $(t)$ of the $x,y,z$ axes, and passing through 
the origin $O$ of the axes. The plane $\Pi$, shown in Fig. 1, intersects
the $xOy$ plane along the line $z=0, x+y=0$, it intersect the 
$yOz$ plane along the line $x=0, y+z=0$, and it intersects the
$xOz$ plane along the line $y=0, x+z=0$.
The condition (\ref{9}) is equivalent to
$(x-y)^2+(x-z)^2+(y-z)^2=0$, which for real $x,y,z$ means that 
$x=y=z$, which represents the trisector line $(t)$ of the axes $x,y,z$. 
The trisector line $(t)$ is perpendicular to the plane $\Pi$. Because of
conditions (\ref{8}) and (\ref{9}), the trisector line $(t)$ and the plane
$\Pi$ will be also called nodal line and respectively nodal plane.

It can be shown that 
if $uu^\prime=0$ then either $u=0$, or $u^\prime=0$, or one of the tricomplex
numbers $u, u^\prime$ belongs to the trisector line $(t)$ and the other belongs
to 
the nodal plane $\Pi$.

\section{Geometric representation of tricomplex numbers}

The tricomplex number $x+hy+kz$ can be represented by the point
$P$ of coordinates $(x,y,z)$. If $O$ is the origin of the axes, 
then the projection $s=OQ$ of the line $OP$ on the trisector line $x=y=z$,
which has the unit tangent $(1/\sqrt{3},1/\sqrt{3},1/\sqrt{3})$, is
\begin{equation}
s=\frac{1}{\sqrt{3}}(x+y+z) . 
\label{10}
\end{equation}
The distance $D=PQ$ from $P$ to the trisector line $x=y=z$, calculated as the 
distance from the point $(x,y,z)$ to the point $Q$ of coordinates 
$[(x+y+z)/3,(x+y+z)/3,(x+y+z)/3]$, is
\begin{equation}
D^2=\frac{2}{3}(x^2+y^2+z^2-xy-xz-yz) .
\label{11}
\end{equation}
The quantities $s$ and $D$ are shown in Fig. 2, where the plane through the
point $P$ and perpendicular to the trisector line $(t)$ intersects the $x$ axis
at point $A$ of coordinates $(x+y+z,0,0)$, the $y$ axis at point $B$ of
coordinates $(0,x+y+z,0)$, and the $z$ axis at point $C$ of coordinates
$(0,0,x+y+z)$.  The azimuthal angle $\phi$ of the tricomplex number $x+hy+kz$
is defined as the angle in the plane $\Pi$ of the projection of $P$ on this
plane, measured from the line of intersection of the plane determined by the
line $(t)$ and the x axis with the plane $\Pi$, $0\leq\phi<2\pi$. 
The expression of $\phi$ in terms of $x,y,z$
can be obtained in a system of coordinates defined by the unit vectors
\begin{equation}
\xi_1: \frac{1}{\sqrt{6}}(2,-1,-1) ; \: \xi_2: \frac{1}{\sqrt{2}}(0,1,-1) ;
\: \xi_3: \frac{1}{\sqrt{3}}(1,1,1) ,
\label{12}
\end{equation}
and having the point $O$ as origin.
The relation between the
coordinates of $P$ in the systems $\xi_1,\xi_2,\xi_3$ and $x,y,z$ can be
written in the form
\begin{equation}
\left(\begin{array}{c}
\xi_1\\ \xi_2\\ \xi_3\end{array}\right)= 
\left(\begin{array}{ccc}
\frac{2}{\sqrt{6}}&-\frac{1}{\sqrt{6}}&-\frac{1}{\sqrt{6}}\\
0&\frac{1}{\sqrt{2}}&-\frac{1}{\sqrt{2}}\\
\frac{1}{\sqrt{3}}&\frac{1}{\sqrt{3}}&\frac{1}{\sqrt{3}}
\end{array}\right)
\left(\begin{array}{c}
x\\
y\\
z
\end{array}\right) .
\label{13}
\end{equation}
The components of the vector $OP$ in the system $\xi_1,\xi_2,\xi_3$ 
can be obtained with the aid of Eq. (\ref{13}) as  
\begin{equation}
(\xi_1,\xi_2,\xi_3)=\left(\frac{1}{\sqrt{6}}(2x-y-z),
\frac{1}{\sqrt{2}}(y-z),\frac{1}{\sqrt{3}}(x+y+z)\right) .
\label{14}
\end{equation}
The expression of the angle $\phi$ as a function of $x,y,z$ is then
\begin{equation}
\cos\phi=\frac{2x-y-z}{2(x^2+y^2+z^2-xy-xz-yz)^{1/2}} ,
\label{15}
\end{equation}
\begin{equation}
\sin\phi=\frac{\sqrt{3}(y-z)}{2(x^2+y^2+z^2-xy-xz-yz)^{1/2}} .
\label{16}
\end{equation}
It can be seen from Eqs. (\ref{15}),(\ref{16}) that the angle of points on the
$x$ axis is $\phi=0$, the angle of points on the $y$ axis is $\phi=2\pi/3$, and
the angle of points on the $z$ axis is $\phi=4\pi/3$.
The angle
$\phi$ is shown in 
Fig. 2 in the plane parallel to $\Pi$, passing through $P$. The axis
$Q\xi_1^{\parallel}$ is parallel to the axis $O\xi_1$, 
the axis $Q\xi_2^{\parallel}$ is parallel to the axis $O\xi_2$, 
and the axis $Q\xi_3^{\parallel}$ is parallel to the axis $O\xi_3$, 
so that, in the plane $ABC$, the angle $\phi$ is measured from the line $QA$.
The angle $\theta$ between the line $OP$ and the trisector line $(t)$ is given
by 
\begin{equation}
\tan\theta=\frac{D}{s},
\label{16a}
\end{equation}
where $0\leq\theta\leq\pi$.
It can be checked that
\begin{equation}
d^2=D^2+s^2 ,
\label{16b}
\end{equation}
where 
\begin{equation}
d^2=x^2+y^2+z^2,
\label{16c}
\end{equation}
so that
\begin{equation}
D=d\sin\theta,\;s=d\cos\theta .
\label{16d}
\end{equation}

The relations (\ref{10}), (\ref{11}), (\ref{15})-(\ref{16a}) can be used to
determine the associated 
projection $s$, the distance $D$, the polar angle $\theta$ with the trisector
line 
$(t)$  and the angle $\phi$ in the
$\Pi$ plane for the tricomplex number $x+hy+kz$.
It can be shown that if $u_1=x_1+hy_1+kz_1, 
u_2=x_2+hy_2+kz_2$ are tricomplex
numbers of projections, distances and angles $s_1,D_1,\theta_1, \phi_1$ and
respectively 
$s_2, D_2, \theta_2, \phi_2$, then the projection $s$, distance $D$ and the
angle $\theta, \phi$ for 
the product tricomplex number $u_1u_2
=x_1x_2+y_1z_2+y_2z_1+h(z_1z_2+x_1y_2+y_1x_2)+k(y_1y_2+x_1z_2+z_1x_2)$
are 
\begin{equation}
s=\sqrt{3}s_1s_2, \: D=\sqrt{\frac{3}{2}}D_1D_2, 
\: \tan\theta=\frac{1}{\sqrt{2}}\tan\theta_1\tan\theta_2\:, \phi=\phi_1+\phi_2.
\label{17}
\end{equation}
The relations (\ref{17}) are consequences of the identities
\begin{eqnarray}
\lefteqn{(x_1x_2+y_1z_2+y_2z_1)+(z_1z_2+x_1y_2+y_1x_2)+(y_1y_2+x_1z_2+z_1x_2)}
\nonumber\\
 & & =(x_1+y_1+z_1)(x_2+y_2+z_2) , 
\label{18}
\end{eqnarray}
\begin{eqnarray}
\lefteqn{(x_1x_2+y_1z_2+y_2z_1)^2+(z_1z_2+x_1y_2+y_1x_2)^2
+(y_1y_2+x_1z_2+z_1x_2)^2}
\nonumber\\
 & & -(x_1x_2+y_1z_2+y_2z_1)(z_1z_2+x_1y_2+y_1x_2)
-(x_1x_2+y_1z_2+y_2z_1)(y_1y_2+x_1z_2+z_1x_2) \nonumber\\
 & & -(z_1z_2+x_1y_2+y_1x_2)+(y_1y_2+x_1z_2+z_1x_2)  \nonumber\\
 & & =(x_1^2+y_1^2+z_1^2-x_1y_1-x_1z_1-y_1z_1) 
 (x_2^2+y_2^2+z_2^2-x_2y_2-x_2z_2-y_2z_2) ,
\label{19}
\end{eqnarray}
\begin{eqnarray}
\lefteqn{\frac{2x_1-y_1-z_1}{2}\frac{2x_2-y_2-z_2}{2}
-\frac{\sqrt{3}}{2}(y_1-z_1)
\frac{\sqrt{3}}{2}(y_2-z_2)}\nonumber\\
& & =\frac{1}{2}[2(x_1x_2+y_1z_2+z_1y_2)
-(z_1z_2+x_1y_2+y_1x_2)-(y_1y_2+x_1z_2+z_1x_2)] ,
\label{20}
\end{eqnarray}
\begin{eqnarray}
\lefteqn{\frac{\sqrt{3}}{2}(y_1-z_1)\frac{2x_2-y_2-z_2}{2}
+\frac{\sqrt{3}}{2}(y_2-z_2)\frac{2x_1-y_1-z_1}{2}}\nonumber\\
& & =\frac{\sqrt{3}}{2} [(z_1z_2+x_1y_2+y_1x_2)-(y_1y_2+x_1z_2+z_1x_2)] .
\label{21}
\end{eqnarray}
The relation (\ref{18}) shows that if $u$ is in the plane
$\Pi$, such that $x+y+z=0$, 
then the product $uu^\prime$ is also in the plane $\Pi$ for any $u^\prime$.
The relation (\ref{19}) shows that if $u$
is on the trisector line $(t)$, such that $x^2+y^2+z^2-xy-xz-yz=0$,
then $uu^\prime$ is also on the trisector line $(t)$ for any $u^\prime$.
If $u, u^\prime$ are points in the plane $x+y+z=1$, then the product
$uu^\prime$ is also in that plane, and if $u, u^\prime$ are points of the
cylindrical surface $x^2+y^2+z^2-xy-xz-yz=1$, then $uu^\prime$ is also in that
cylindrical surface. This means that if $u, u^\prime$ are points on the circle
$x+y+z=1, x^2+y^2+z^2-xy-xz-yz=1,$  
which is perpendicular to the trisector line, is situated at
a distance $1/\sqrt{3}$ from the origin and has the radius $\sqrt{2/3}$,
then the tricomplex product $uu^\prime$ is also on the same circle. This
invariant circle for the multiplication of tricomplex numbers is described by
the equations  
\begin{equation}
x=\frac{1}{3}+\frac{2}{3}\cos\phi,\:
y=\frac{1}{3}-\frac{1}{3}\cos\phi+\frac{1}{\sqrt{3}}\sin\phi,\:
z=\frac{1}{3}-\frac{1}{3}\cos\phi-\frac{1}{\sqrt{3}}\sin\phi.
\label{21b}
\end{equation}
It has the center at the point (1/3,1/3,1/3) and passes through the points
(1,0,0), (0,1,0) and (0,0,1), as shown in Fig. 3. 

An important quantity is the amplitude $\rho$ defined as
$\rho=\nu^{1/3}$, so that
\begin{equation}
\rho^3=x^3+y^3+z^3-3xyz .
\label{22}
\end{equation}
The amplitude $\rho$ of the product $u_1u_2$ of the tricomplex numbers
$u_1,u_2$ of amplitudes $\rho_1,\rho_2$ is
\begin{equation}
\rho=\rho_1\rho_2 ,
\label{23}
\end{equation}
as can be seen from the identity
\begin{eqnarray}
\lefteqn{(x_1x_2+y_1z_2+y_2z_1)^3+(z_1z_2+x_1y_2+y_1x_2)^3
+(y_1y_2+x_1z_2+z_1x_2)^3} 
\nonumber\\
 &&-3(x_1x_2+y_1z_2+y_2z_1)(z_1z_2+x_1y_2+y_1x_2)(y_1y_2+x_1z_2+z_1x_2)
\nonumber\\ 
 &&=(x_1^3+y_1^3+z_1^3-3x_1y_1z_1)(x_2^3+y_2^3+z_2^3-3x_2y_2z_2) .
\label{24}
\end{eqnarray}
The identity in Eq. (\ref{24}) can be demonstrated with the aid of 
Eqs. (\ref{7}), (\ref{18}) and (\ref{19}). Another method
would be to use the representation of the multiplication of the tricomplex
numbers by matrices, in which the tricomplex number $u=x+hy+kz$ is
represented by the matrix
\begin{equation}
\left(\begin{array}{ccc}
x&y&z\\
z&x&y\\
y&z&x 
\end{array}\right) .
\label{25}
\end{equation}
The product $u=x+hy+kz$ of the tricomplex numbers $u_1=x_1+hy_1+kz_1,
u_2=x_2+hy_2+kz_2$, is represented by the matrix multiplication
\begin{equation}
\left(\begin{array}{ccc}
x&y&z\\
z&x&y\\
y&z&x 
\end{array}\right) =
\left(\begin{array}{ccc}
x_1&y_1&z_1\\
z_1&x_1&y_1\\
y_1&z_1&x_1 
\end{array}\right) 
\left(\begin{array}{ccc}
x_2&y_2&z_2\\
z_2&x_2&y_2\\
y_2&z_2&x_2 
\end{array}\right) .
\label{26}
\end{equation}
If 
\begin{equation}
\nu={\rm det}\left(\begin{array}{ccc}
x&y&z\\
z&x&y\\
y&z&x 
\end{array}\right),
\label{27}
\end{equation}
it can be checked that
\begin{equation}
\nu=x^3+y^3+z^3-3xyz .
\label{27b}
\end{equation}
The identity (\ref{24}) is then a consequence of the fact the determinant 
of the product of matrices is equal to the product of the determinants 
of the factor matrices. 

It can be seen from Eqs. (\ref{10}) and (\ref{11}) that
\begin{equation}
x^3+y^3+z^3-3xyz=\frac{3\sqrt{3}}{2}sD^2 ,
\label{28}
\end{equation}
which can be written with the aid of relations (\ref{16d}) and (\ref{22}) as
\begin{equation}
\rho=\frac{3^{1/2}}{2^{1/3}}d\sin^{2/3}\theta\cos^{1/3}\theta.
\label{28b}
\end{equation}
This means that the surfaces of constant $\rho$ are surfaces of rotation
having the trisector line $(t)$ as axis, as shown in Fig. 4.

\section{The tricomplex cosexponential functions}

The exponential function of the tricomplex variable $u$ can be defined by the
series
\begin{equation}
\exp u = 1+u+u^2/2!+u^3/3!+\cdots . 
\label{29}
\end{equation}
It can be checked by direct multiplication of the series that
\begin{equation}
\exp(u+u^\prime)=\exp u \cdot \exp u^\prime ,
\label{30}
\end{equation}
which is valid as long as the multiplication is a commutative operation.
If $u=x+hy+kz$, then $\exp u$ can be calculated as $\exp u=\exp x \cdot
\exp (hy) \cdot \exp (kz)$. According to Eq. (\ref{1}), $h^2=k, h^3=1, k^2=h,
k^3=1$, and in general 
\begin{equation}
h^{3m}=1, h^{3m+1}=h, h^{3m+2}=k, k^{3m}=1, k^{3m+1}=k, k^{3m+2}=h ,
\label{31}
\end{equation}
where $n$ is a natural number,
so that $\exp (hy)$ and $\exp(kz)$ can be written as
\begin{equation}
\exp (hy) = {\rm cx} \:y + h\: {\rm mx} \:y + k\: {\rm px} \:y , 
\label{32}
\end{equation}
\begin{equation}
\exp (kz) = {\rm cx} \:z + h \:{\rm px} \:z + k \:{\rm mx} \:z , 
\label{33}
\end{equation}
where the functions cx, mx, px, which will be called in this work polar
cosexponential functions, are defined by the series 
\begin{equation}
{\rm cx \:y} = 1 + y^3/3! +y^6/6!+\cdots
\label{34}
\end{equation}
\begin{equation}
{\rm mx} \:y = y + y^4/4! + y^7/7!+\cdots
\label{35}
\end{equation}
\begin{equation}
{\rm px} \:y = y^2/2! + y^5/5! +y^8/8!+\cdots .
\label{36}
\end{equation}
From the series definitions it can be seen that ${\rm cx} \: 0 =1, 
{\rm mx} \: 0=0, {\rm px} \: 0=0.$
The tridimensional polar cosexponential functions belong to the class
of the polar n-dimensional cosexponential functions
$g_{nk}$, \cite{2d} and ${\rm cx}=g_{30}, {\rm mx}=g_{31}, {\rm px}=g_{32}$.
It can be checked that
\begin{equation}
{\rm cx} \: y +{\rm px} \: y+{\rm mx} \: y =\exp y .
\label{36a}
\end{equation}
By expressing the fact that $\exp(hy+hz)=\exp(hy)\cdot\exp(hz)$ with the aid of
the cosexponential functions (\ref{34})-(\ref{36}) the following addition
theorems can be obtained
\begin{equation}
{\rm cx}\:(y+z)={\rm cx}\:y \:{\rm cx}\:z+{\rm mx}\:y 
\:{\rm px}\:z + {\rm px}\:y \:{\rm mx}\:z  ,
\label{38}
\end{equation}
\begin{equation}
{\rm mx}\:(y+z)={\rm px}\:y \:{\rm px}\:z+{\rm cx}\:y 
\:{\rm mx}\:z + {\rm mx}\:y \:{\rm cx}\:z  ,
\label{39}
\end{equation}
\begin{equation}
{\rm px}\:(y+z)={\rm mx}\:y \:{\rm mx}\:z+{\rm cx}\:y 
\:{\rm px}\:z + {\rm px}\:y \:{\rm cx}\:z  .
\label{40}
\end{equation}
For $y=z$, Eqs. (\ref{38})-(\ref{40}) yield
\begin{equation}
{\rm cx}\: 2y={\rm cx}^2\:y +2 \:{\rm mx}\:y 
\:{\rm px}\:z ,
\label{38a}
\end{equation}
\begin{equation}
{\rm mx}\:2y={\rm px}^2 \:y +2\:{\rm cx}\:y 
\:{\rm mx}\:z ,
\label{39a}
\end{equation}
\begin{equation}
{\rm px}\:2y={\rm mx}^2\:y +2\:{\rm cx}\:y 
\:{\rm px}\:z .
\label{40a}
\end{equation}
The cosexponential functions are neither even nor odd functions. 
For $z=-y$, Eqs. (\ref{38})-(\ref{40}) yield
\begin{equation}
{\rm cx}\:y \:{\rm cx}\:(-y)+{\rm mx}\:y 
\:{\rm px}\:(-y) + {\rm px}\:y \:{\rm mx}\:(-y) =1  ,
\label{38b}
\end{equation}
\begin{equation}
{\rm px}\:y \:{\rm px}\:(-y)+{\rm cx}\:y 
\:{\rm mx}\:(-y) + {\rm mx}\:y \:{\rm cx}\:(-y)=0  ,
\label{39b}
\end{equation}
\begin{equation}
{\rm mx}\:y \:{\rm mx}\:(-y)+{\rm cx}\:y 
\:{\rm px}\:(-y) + {\rm px}\:y \:{\rm cx}\:(-y)=0  .
\label{40b}
\end{equation}

Expressions of the cosexponential functions in terms of regular exponential and
cosine functions can be obtained by considering the series expansions for
$e^{(h+k)y}$ and $e^{(h-k)y}$. These expressions can be obtained by calculating
first $(h+k)^n$ and $(h-k)^n$. It can be shown that
\begin{equation}
(h+k)^m=\frac{1}{3}\left[(-1)^{m-1}+2^m\right](h+k)+
\frac{2}{3}\left[(-1)^m+2^{m-1}\right] ,
\label{41}
\end{equation}
\begin{equation}
(h-k)^{2m}=(-1)^{m-1}3^{m-1}(k+k-2), \:(h-k)^{2m+1}=(-1)^m 3^m (h-k) ,
\label{42}
\end{equation}
where $n$ is a natural number. Then
\begin{equation}
e^{(h+k)y}=(h+k)\left(-\frac{1}{3}e^{-y}+\frac{1}{3}e^{2y}\right)
+\frac{2}{3}e^{-y}+\frac{1}{3}e^{2y} .
\label{43}
\end{equation}
As a corollary, the following identities can be obtained from Eq. (\ref{43}) by
writing $e^{(h+k)y}=e^{hy}e^{ky}$ and expressing $e^{hy}$ and $e^{ky}$ in terms
of cosexponential functions via Eqs. (\ref{32}) and (\ref{33}),
\begin{equation}
{\rm cx}^2 \: y +{\rm mx}^2 \: y+{\rm px}^2 \: y=\frac{2}{3} e^{-y}
+\frac{1}{3}e^{2y} , 
\label{43a}
\end{equation}
\begin{equation}
 {\rm cx} \: y \:\: {\rm mx} \: y
+ {\rm cx} \: y \:\: {\rm px} \: y
+ {\rm mx} \: y \:\: {\rm px} \: y  = -\frac{1}{3}e^{-y}+\frac{1}{3}
e^{2y} . 
\label{43b}
\end{equation}
From Eqs. (\ref{43a}) and (\ref{43b}) it results that
\begin{eqnarray}
\lefteqn{{\rm cx}^2 \: y +{\rm mx}^2 \: y+{\rm px}^2 \: y
}\nonumber\\ 
& & - {\rm cx} \: y \:\: {\rm mx} \: y
- {\rm cx} \: y \:\: {\rm px} \: y
- {\rm mx} \: y \:\: {\rm px} \: y  = \exp(-y) . 
\label{36b}
\end{eqnarray}
Then from Eqs. (\ref{7}), (\ref{36a}) and (\ref{36b}) it follows that
\begin{equation}
{\rm cx}^3 \: y +{\rm mx}^3 \: y+{\rm px}^3 \: y
 - 3 {\rm cx} \: y \:\: {\rm mx} \: y \:\:{\rm px} \: y =1 .
\label{37}
\end{equation}

Similarly,
\begin{equation}
e^{(h-k)y}=\frac{1}{3}(1+h+k)+\frac{1}{3}(2-h-k)\cos(\sqrt{3}y)
+\frac{1}{\sqrt{3}}(h-k)\sin(\sqrt{3}y) .
\label{44}
\end{equation}
The last relation can also be written as
\begin{eqnarray}
e^{(h-k)y}&=&\frac{1}{3}+\frac{2}{3}\cos(\sqrt{3}y)+h\left[\frac{1}{3}
+\frac{2}{3}\cos\left(\sqrt{3}y-\frac{2\pi}{3}\right)\right]\nonumber\\
&&+k\left[\frac{1}{3}+\frac{2}{3}\cos\left(\sqrt{3}y+\frac{2\pi}{3}
\right)\right].
\label{44a}
\end{eqnarray}
As a corollary, the following identities can be obtained from Eq. (\ref{44}) by
writing $e^{(h-k)y}=e^{hy}e^{-ky}$ and expressing $e^{hy}$ and $e^{-ky}$ in
terms 
of cosexponential functions via Eqs. (\ref{32}) and (\ref{33}),
\begin{equation}
 {\rm cx} \: y \:\: {\rm cx} \: (-y)
+ {\rm mx} \: y \:\: {\rm mx} \: (-y)
+ {\rm px} \: y \:\: {\rm px} \: (-y)  =
\frac{1}{3}+\frac{2}{3}\cos(\sqrt{3}y) , 
\label{44b}
\end{equation}
\begin{eqnarray}
 {\rm cx} \: y \:\: {\rm px} \: (-y)
+ {\rm mx} \: y \:\: {\rm cx} \: (-y)
+ {\rm px} \: y \:\: {\rm mx} \: (-y) 
=\frac{1}{3}+\frac{2}{3}\cos\left(\sqrt{3}y-\frac{2\pi}{3}\right)
\label{44d}
\end{eqnarray}
\begin{eqnarray}
 {\rm cx} \: y \:\: {\rm mx} \: (-y)
+ {\rm mx} \: y \:\: {\rm px} \: (-y)
+ {\rm px} \: y \:\: {\rm cx} \: (-y) 
=\frac{1}{3}+\frac{2}{3}\cos\left(\sqrt{3}y+\frac{2\pi}{3}\right)
\label{44c}
\end{eqnarray}

Expressions of $e^{2hy}$ in terms of the regular exponential and
cosine functions can be obtained by the multiplication of the expressions of
$e^{(h+k)y}$ and $e^{(h-k)y}$ from Eqs. (\ref{43}) and (\ref{44}).
At the same time, Eq. (\ref{32}) gives an expression of $e^{2hy}$ in terms 
of cosexponential functions. By equating the real and hypercomplex parts of 
these two forms of $e^{2y}$ and then
replacing $2y$ by $y$ gives the expressions of the cosexponential functions as
\begin{equation}
{\rm cx} \:y =\frac{1}{3}\: e^y+\frac{2}{3}\cos\left(\frac{\sqrt{3}}{2}y\right)
\:e^{-y/2},
\label{45}
\end{equation}
\begin{equation}
{\rm mx} \:y =\frac{1}{3}\: e^y+\frac{2}{3}\cos\left(\frac{\sqrt{3}}{2}y
-\frac{2\pi}{3}\right)
\:e^{-y/2} ,
\label{46}
\end{equation}
\begin{equation}
{\rm px} \:y =\frac{1}{3}\: e^y+\frac{2}{3}\cos\left(\frac{\sqrt{3}}{2}y
+\frac{2\pi}{3}\right)
\:e^{-y/2} .
 \label{47}
\end{equation}
It is remarkable that the series in Eqs. (\ref{34})-(\ref{36}), in which the 
terms are either of the form $y^{3m}$, or $y^{3m+1}$, or $y^{3m+2}$, 
can be expressed
in terms of elementary functions whose power series are not subject to such 
restrictions. 
The cosexponential functions differ by the
phase of the cosine function in their expression, and the designation of the
functions in Eqs. (\ref{46}) and (\ref{47}) as mx and px
refers respectively to the minus or plus sign of the phase term $2\pi/3$. 
The graphs of the cosexponential functions are shown in Fig. 5.

It can be checked that the cosexponential functions are solutions of the
third-order differential equation
\begin{equation}
\frac{d^3\zeta}{du^3}=\zeta ,
\label{47b}
\end{equation}
whose solutions are of the form $\zeta(u)=A\:{\rm cx}\: u+B\:{\rm mx}\: u
+C\:{\rm px}\: u.$
It can also be checked that the derivatives of the cosexponential functions are
related by
\begin{equation}
\frac{d{\rm px}}{du}={\rm mx}, \:
\frac{d{\rm mx}}{du}={\rm cx}, \:
\frac{d{\rm cx}}{du}={\rm px} .
\label{47c}
\end{equation}

\section{Exponential and trigonometric forms of tricomplex numbers}

If for a tricomplex number $u=x+ky+kz$ another tricomplex number 
$u_1=x_1+hy_1+kz_1$ exists 
such that
\begin{equation}
x+hy+kz=e^{x_1+hy_1+kz_1} ,
\label{48}
\end{equation}
then $u_1$ is said to be the logarithm of $u$,
\begin{equation}
u_1=\ln u .
\label{49}
\end{equation}
The expressions of $x_1, y_1, z_1$ as functions of $x, y, z$ can be obtained by
developing $e^{hy_1}$ and $e^{kz_1}$ with the aid of Eqs. (\ref{32}) and
(\ref{33}), by multiplying these expressions and separating the hypercomplex
components, 
\begin{equation}
x=e^{x_1}[ {\rm cx} \: y_1 \:\: {\rm cx} \: z_1
+ {\rm mx} \: y_1 \:\: {\rm mx} \: z_1
+ {\rm px} \: y_1 \:\: {\rm px} \: z_1 ] ,
\label{50}
\end{equation}
\begin{equation}
y=e^{x_1}[ {\rm cx} \: y_1 \:\: {\rm px} \: z_1
+ {\rm mx} \: y_1 \:\: {\rm cx} \: z_1
+ {\rm px} \: y_1 \:\: {\rm mx} \: z_1 ] ,
\label{51}
\end{equation}
\begin{equation}
z=e^{x_1}[ {\rm cx} \: y_1 \:\: {\rm mx} \: z_1
+ {\rm px} \: y_1 \:\: {\rm cx} \: z_1
+ {\rm mx} \: y_1 \:\: {\rm px} \: z_1 ] ,
\label{52}
\end{equation}
Using Eq. (\ref{24}) with the substitutions $x_1\rightarrow{\rm cx}\: y_1,
y_1\rightarrow{\rm mx}\: y_1, z_1\rightarrow{\rm px}\: y_1, 
x_2\rightarrow{\rm cx}\: z_1, y_2\rightarrow{\rm px}\: z_1,
z_2\rightarrow{\rm mx}\: z_1$ and then the identity (\ref{37}) yields
\begin{equation}
x^3+y^3+z^3-3xyz=e^{3x_1} , 
\label{53}
\end{equation}
whence
\begin{equation}
x_1=\frac{1}{3}\ln(x^3+y^3+z^3-3xyz) .
\label{54}
\end{equation}
The logarithm in Eq. (\ref{54}) exists as a real function for $x+y+z>0$.
A further relation can be obtained by summing Eqs. (\ref{50})-(\ref{52}) and
then using the addition theorems (\ref{38})-(\ref{40}) 
\begin{equation}
\frac{x+y+z}{(x^3+y^3+z^3-3xyz)^{1/3}}={\rm cx}\:(y_1+z_1)+
{\rm mx}\:(y_1+z_1)+{\rm px}\:(y_1+z_1) .
\label{55}
\end{equation}
The sum in Eq. (\ref{55}) is according to Eq. (\ref{36a}) $e^{y_1+z_1}$, so
that 
\begin{equation}
y_1+z_1=\ln\frac{x+y+z}{(x^3+y^3+z^3-3xyz)^{1/3}} .
\label{56}
\end{equation}
The logarithm in Eq. (\ref{56}) is defined for points which are not on the
trisector line $(t)$, so that $x^2+y^2+z^2-xy-xz-yz\not= 0$.
Substituting in Eq. (\ref{48}) the expression of $x_1$, Eq. (\ref{54}), and of
$z_1$ as a function of $x, y, z, y_1,$ Eq. (\ref{56}), yields
\begin{equation}
\frac{u}{\rho}
\exp\left[-k\ln\left(\frac{\sqrt{2}s}{D}\right)^{2/3}\right]=e^{(h-k)y_1} ,
\label{57}
\end{equation}
where the quantities $\rho, s$ and $D$ have been defined in Eqs.
(\ref{22}),(\ref{10}) and (\ref{11}).
Developing the exponential functions in the left-hand side of Eq. (\ref{57}) 
with the aid of Eq. (\ref{33}) and using the expressions of the cosexponential
functions, Eqs. (\ref{45})-(\ref{47}), and using the relation (\ref{44}) for
the right-hand side of Eq. (\ref{57}) yields for the real part
\begin{equation}
\frac{\left(x-\frac{y+z}{2}\right)\cos\left[\frac{1}{\sqrt{3}}\ln
\left(\frac{\sqrt{2}s}{D}\right)\right]
-\frac{\sqrt{3}}{2}(y-z)\sin
\left[\frac{1}{\sqrt{3}}\ln
\left(\frac{\sqrt{2}s}{D}\right)\right]}
{(x^2+y^2+z^2-xy-xz-yz)^{1/2}}=\cos(\sqrt{3}y_1) ,
\label{58}
\end{equation}
which can also be written as
\begin{equation}
\cos\left[\frac{1}{\sqrt{3}}\ln\left(\frac{\sqrt{2}s}{D}\right)+\phi\right]
=\cos(\sqrt{3}y_1)
\label{59}
\end{equation}
where $\phi$ is the angle defined in Eqs. (\ref{15}) and (\ref{16}).
Thus
\begin{equation}
y_1=\frac{1}{3}\ln\left(\frac{\sqrt{2}s}{D}\right)+\frac{1}{\sqrt{3}}\phi .
\label{60}
\end{equation}
The exponential form of the tricomplex number $u$ is then
\begin{equation}
u=\rho\:
\exp\left[\frac{1}{3}(h+k)\ln\frac{\sqrt{2}}{\tan\theta}+
\frac{1}{\sqrt{3}}(h-k)\phi\right] ,
\label{61}
\end{equation}
where $\theta$ is the angle between the line $OP$ 
connecting the origin to the point $P$ of coordinates $(x,y,z)$ and the
trisector line $(t)$, defined in Eq. (\ref{16a}) and shown in Fig. 2.
The exponential in Eq. (\ref{61}) can be expanded with the aid of Eq.
(\ref{43}) and (\ref{44a}) as
\begin{equation}
\exp\left[\frac{1}{3}(h+k)\ln\frac{\sqrt{2}}{\tan\theta}\right] 
=\frac{2-h-k}{3}\left(\frac{\tan\theta}{\sqrt{2}}\right)^{1/3}
+\frac{1+h+k}{3}\left(\frac{\sqrt{2}}{\tan\theta}\right)^{2/3} ,
\label{62}
\end{equation}
so that 
\begin{equation}
x+hy+kz=\rho
\left[\frac{2-h-k}{3}\left(\frac{\tan\theta}{\sqrt{2}}\right)^{1/3}
+\frac{1+h+k}{3}\left(\frac{\sqrt{2}}{\tan\theta}\right)^{2/3}\right]
\exp\left(\frac{h-k}{\sqrt{3}}\phi\right) .
\label{63}
\end{equation}
Substituting in Eq. (\ref{63}) the expression of the amplitude $\rho$, Eq.
(\ref{28b}), yields
\begin{equation}
u=d\sqrt{\frac{3}{2}}
\left(\frac{2-h-k}{3}\sin\theta
+\frac{1+h+k}{3}\sqrt{2}\cos\theta\right)
\exp\left(\frac{h-k}{\sqrt{3}}\phi\right) ,
\label{63b}
\end{equation}
which is the trigonometric form of the tricomplex number $u$.
As can be seen from Eq. (\ref{63b}), the tricomplex number $x+hy+kz$ is
written as the product of the modulus $d$, of a factor depending on the
polar angle $\theta$ with respect to the trisector line, and of a factor 
depending 
of the azimuthal angle $\phi$ in the plane $\Pi$ perpendicular to the
trisector line. 
The exponential in Eq. (\ref{63b}) can be expanded further with the aid of
Eq. (\ref{44a}) as 
\begin{equation}
\exp\left(\frac{1}{\sqrt{3}}(h-k)\phi\right) 
=\frac{1+h+k}{3}+\frac{2-h-k}{3}\cos\phi +\frac{h-k}{\sqrt{3}}\sin\phi ,
\label{64}
\end{equation}
so that the tricomplex number $x+hy+kz$ can also be written, after
multiplication of the factors, in the form
\begin{eqnarray}
\lefteqn{x+hy+kz=
\frac{2-h-k}{3}(x^2+y^2+z^2-xy-xz-yz)^{1/2}\cos\phi}\nonumber\\
&&+\frac{h-k}{\sqrt{3}}(x^2+y^2+z^2-xy-xz-yz)^{1/2}\sin\phi
+\frac{1+h+k}{3}(x+y+z)
\label{65}
\end{eqnarray}
The validity of Eq. (\ref{65}) can be checked by substituting the expressions
of $\cos\phi$ and $\sin\phi$ from Eqs. (\ref{15}) and (\ref{16}).
 
\section{Elementary functions of a tricomplex variable}

It can be shown with the aid of Eq. (\ref{61}) that 
\begin{equation}
(x+hy+kz)^m=
\rho^m
\left[\frac{2-h-k}{3}\left(\frac{\tan\theta}{\sqrt{2}}\right)^{m/3}
+\frac{1+h+k}{3}\left(\frac{\sqrt{2}}{\tan\theta}\right)^{2m/3}\right]
\exp\left(\frac{h-k}{\sqrt{3}}m\phi\right) ,
\label{66}
\end{equation}
or equivalently
\begin{eqnarray}
\lefteqn{(x+hy+kz)^m=
\frac{2-h-k}{3}(x^2+y^2+z^2-xy-xz-yz)^{m/2}\cos(m\phi)}
\nonumber\\
&&+\frac{h-k}{\sqrt{3}}(x^2+y^2+z^2-xy-xz-yz)^{m/2}\sin(m\phi)
+\frac{1+h+k}{3}(x+y+z)^m
\label{67}
\end{eqnarray}
which are valid for real values of $m$. Thus Eqs. (\ref{66}) or (\ref{67})
define the power function $u^m$ of the tricomplex variable
$u$. 

The power function is multivalued unless $m$ is an integer. 
It can be inferred from Eq. (\ref{61}) that, for integer values of $m$,
\begin{equation}
(uu^\prime)^m=u^m\:u^{\prime m} .
\label{68}
\end{equation}
For natural
$m$, Eq. (\ref{67}) can be checked with the aid of relations (\ref{15}) and
(\ref{16}). For example if $m=2$, it can be checked that the right-hand side of
Eq. (\ref{67}) is equal to $(x+hy+kz)^2=x^2+2yz+h(z^2+2xz)+k(y^2+2xz)$.

The logarithm $u_1$ of the tricomplex number $u$, $u_1=\ln u$, can be defined
as the solution of Eq. (\ref{48}) for $u_1$ as a function of $u$. 
For $x+y+z>0$, from Eq. (\ref{61}) it results that 
\begin{equation}
\ln u=\ln \rho+
\frac{1}{3}(h+k)\ln\left(\frac{\tan\theta}{\sqrt{2}}\right)+
\frac{1}{\sqrt{3}}(h-k)\phi .
\label{69}
\end{equation}
It can be checked with the aid of Eqs. (\ref{17}) and (\ref{23}) that
\begin{equation}
\ln(uu^\prime)=\ln u+\ln u^\prime ,
\label{69b}
\end{equation}
which is valid up to integer multiples of $2\pi(h-k)/\sqrt{3}$.

The trigonometric functions $\cos u$ and $\sin u $ of the tricomplex variable
$u$ are defined by the series
\begin{equation}
\cos u = 1 - u^2/2!+u^4/4!+\cdots, 
\label{70}
\end{equation}
\begin{equation}
\sin u=u-u^3/3!+u^5/5! +\cdots .
\label{71}
\end{equation}
It can be checked by series multiplication that the usual addition theorems
hold also for the tricomplex numbers $u, u^\prime$,
\begin{equation}
\cos(u+u^\prime)=\cos u\cos u^\prime - \sin u \sin u^\prime ,
\label{72}
\end{equation}
\begin{equation}
\sin(u+u^\prime)=\sin u\cos u^\prime + \cos u \sin u^\prime .
\label{73}
\end{equation}
The trigonometric functions of the hypercomplex variables $hy, ky$ can be
expressed in terms of the cosexponential functions as
\begin{equation}
\cos(hy)=\frac{1}{2}[{\rm cx}\:(iy)+{\rm cx}\: (-iy)]+ 
\frac{1}{2}h[{\rm mx}\:(iy)+{\rm mx}\: (-iy)]+
\frac{1}{2}k[{\rm px}\:(iy)+{\rm px}\: (-iy)] ,
\label{74}
\end{equation}
\begin{equation}
\cos(ky)=\frac{1}{2}[{\rm cx}\:(iy)+{\rm cx}\: (-iy)]+ 
\frac{1}{2}h[{\rm px}\:(iy)+{\rm px}\: (-iy)]+
\frac{1}{2}k[{\rm mx}\:(iy)+{\rm mx}\: (-iy)] ,
\label{75}
\end{equation}
\begin{equation}
\sin(hy)=\frac{1}{2i}[{\rm cx}\:(iy)-{\rm cx}\: (-iy)]+ 
\frac{1}{2i}h[{\rm mx}\:(iy)-{\rm mx}\: (-iy)]+
\frac{1}{2i}k[{\rm px}\:(iy)-{\rm px}\: (-iy)] ,
\label{76}
\end{equation}
\begin{equation}
\sin(ky)=\frac{1}{2i}[{\rm cx}\:(iy)-{\rm cx}\: (-iy)]+ 
\frac{1}{2i}h[{\rm px}\:(iy)-{\rm px}\: (-iy)]+
\frac{1}{2i}k[{\rm mx}\:(iy)-{\rm mx}\: (-iy)] ,
\label{77}
\end{equation}
where $i$ is the imaginary unit. Using the expressions of the cosexponential
functions in Eqs. (\ref{45})-(\ref{47}) gives expressions of the 
trigonometric functions of $hy, hz$ as
\begin{eqnarray}
\lefteqn{\cos(hy) =\frac{1}{3}\cos y+
\frac{2}{3}\cosh\left(\frac{\sqrt{3}}{2}y\right)\cos\frac{y}{2}+}\nonumber\\ 
&&+h\left[\frac{1}{3}\cos y
-\frac{1}{3}\cosh\left(\frac{\sqrt{3}}{2}y\right)\cos\frac{y}{2} 
+\frac{1}{\sqrt{3}}\sinh\left(\frac{\sqrt{3}}{2}y\right)\sin\frac{y}{2}
\nonumber\right] \\
&&+k\left[\frac{1}{3}\cos y
-\frac{1}{3}\cosh\left(\frac{\sqrt{3}}{2}y\right)\cos\frac{y}{2} 
-\frac{1}{\sqrt{3}}\sinh\left(\frac{\sqrt{3}}{2}y\right)\sin\frac{y}{2}
\right]
\label{78}
\end{eqnarray}
\begin{eqnarray}
\lefteqn{\cos(ky) =\frac{1}{3}\cos y+
\frac{2}{3}\cosh\left(\frac{\sqrt{3}}{2}y\right)\cos\frac{y}{2}+
}\nonumber\\ 
&&+h\left[\frac{1}{3}\cos y
-\frac{1}{3}\cosh\left(\frac{\sqrt{3}}{2}y\right)\cos\frac{y}{2} 
-\frac{1}{\sqrt{3}}\sinh\left(\frac{\sqrt{3}}{2}y\right)\sin\frac{y}{2}
\nonumber\right] \\
&&+k\left[\frac{1}{3}\cos y
-\frac{1}{3}\cosh\left(\frac{\sqrt{3}}{2}y\right)\cos\frac{y}{2} 
+\frac{1}{\sqrt{3}}\sinh\left(\frac{\sqrt{3}}{2}y\right)\sin\frac{y}{2}
\right]
\label{79}
\end{eqnarray}
\begin{eqnarray}
\lefteqn{\sin(hy) =\frac{1}{3}\sin y
-\frac{2}{3}\cosh\left(\frac{\sqrt{3}}{2}y\right)\sin\frac{y}{2}}
\nonumber\\ 
&&
+h\left[\frac{1}{3}\sin y
+\frac{1}{3}\cosh\left(\frac{\sqrt{3}}{2}y\right)\sin\frac{y}{2} 
+\frac{1}{\sqrt{3}}\sinh\left(\frac{\sqrt{3}}{2}y\right)\cos\frac{y}{2}
\nonumber\right] \\
&&
+k\left[\frac{1}{3}\sin y
+\frac{1}{3}\cosh\left(\frac{\sqrt{3}}{2}y\right)\sin\frac{y}{2} 
-\frac{1}{\sqrt{3}}\sinh\left(\frac{\sqrt{3}}{2}y\right)\cos\frac{y}{2}
\right]
\label{80}
\end{eqnarray}
\begin{eqnarray}
\lefteqn{\sin(ky) =\frac{1}{3}\sin y
-\frac{2}{3}\cosh\left(\frac{\sqrt{3}}{2}y\right)\sin\frac{y}{2}}\nonumber\\ 
&&
+h\left[\frac{1}{3}\sin y
+\frac{1}{3}\cosh\left(\frac{\sqrt{3}}{2}y\right)\sin\frac{y}{2} 
-\frac{1}{\sqrt{3}}\sinh\left(\frac{\sqrt{3}}{2}y\right)\cos\frac{y}{2}
\nonumber\right] \\
&&
+k\left[\frac{1}{3}\sin y
+\frac{1}{3}\cosh\left(\frac{\sqrt{3}}{2}y\right)\sin\frac{y}{2} 
+\frac{1}{\sqrt{3}}\sinh\left(\frac{\sqrt{3}}{2}y\right)\cos\frac{y}{2}
\right]
\label{81}
\end{eqnarray}
The trigonometric functions of a tricomplex number $x+hy+kz$ can then be
expressed in terms of elementary functions with the aid of the addition
theorems Eqs. (\ref{72}), (\ref{73}) and of the expressions in  Eqs. 
(\ref{78})-(\ref{81}).

The hyperbolic functions $\cosh u$ and $\sinh u $ of the fourcomplex variable
$u$ are defined by the series
\begin{equation}
\cosh u = 1 + u^2/2!+u^4/4!+\cdots, 
\label{81a}
\end{equation}
\begin{equation}
\sinh u=u+u^3/3!+u^5/5! +\cdots .
\label{81b}
\end{equation}
It can be checked by series multiplication that the usual addition theorems
hold also for the fourcomplex numbers $u, u^\prime$,
\begin{equation}
\cosh(u+u^\prime)=\cosh u\cosh u^\prime + \sinh u \sinh u^\prime ,
\label{81c}
\end{equation}
\begin{equation}
\sinh(u+u^\prime)=\sinh u\cosh u^\prime + \cosh u \sinh u^\prime .
\label{81d}
\end{equation}
The hyperbolic functions of the hypercomplex variables $hy, ky$ can be
expressed in terms of the elementary functions as
\begin{eqnarray}
\lefteqn{\cosh(hy) =\frac{1}{3}\cosh y+
\frac{2}{3}\cos\left(\frac{\sqrt{3}}{2}y\right)\cosh\frac{y}{2}+}\nonumber\\ 
&&+h\left[\frac{1}{3}\cosh y
-\frac{1}{3}\cos\left(\frac{\sqrt{3}}{2}y\right)\cosh\frac{y}{2} 
-\frac{1}{\sqrt{3}}\sin\left(\frac{\sqrt{3}}{2}y\right)\sinh\frac{y}{2}
\nonumber\right] \\
&&+k\left[\frac{1}{3}\cosh y
-\frac{1}{3}\cos\left(\frac{\sqrt{3}}{2}y\right)\cosh\frac{y}{2} 
+\frac{1}{\sqrt{3}}\sin\left(\frac{\sqrt{3}}{2}y\right)\sinh\frac{y}{2}\right]
\label{81e}
\end{eqnarray}
\begin{eqnarray}
\lefteqn{\cosh(ky) =\frac{1}{3}\cosh y+
\frac{2}{3}\cos\left(\frac{\sqrt{3}}{2}y\right)\cosh\frac{y}{2}+}\nonumber\\ 
&&+h\left[\frac{1}{3}\cosh y
-\frac{1}{3}\cos\left(\frac{\sqrt{3}}{2}y\right)\cosh\frac{y}{2} 
+\frac{1}{\sqrt{3}}\sin\left(\frac{\sqrt{3}}{2}y\right)\sinh\frac{y}{2}
\nonumber\right] \\
&&+k\left[\frac{1}{3}\cosh y
-\frac{1}{3}\cos\left(\frac{\sqrt{3}}{2}y\right)\cosh\frac{y}{2} 
-\frac{1}{\sqrt{3}}\sin\left(\frac{\sqrt{3}}{2}y\right)\sinh\frac{y}{2}
\right]
\label{81f}
\end{eqnarray}
\begin{eqnarray}
\lefteqn{\sinh(hy) =\frac{1}{3}\sinh y
-\frac{2}{3}\cos\left(\frac{\sqrt{3}}{2}y\right)\sinh\frac{y}{2}
}\nonumber\\ 
&&+h\left[\frac{1}{3}\sinh y
+\frac{1}{3}\cos\left(\frac{\sqrt{3}}{2}y\right)\sinh\frac{y}{2} 
+\frac{1}{\sqrt{3}}\sin\left(\frac{\sqrt{3}}{2}y\right)\cosh\frac{y}{2}
\nonumber\right] \\
&&+k\left[\frac{1}{3}\sinh y
+\frac{1}{3}\cos\left(\frac{\sqrt{3}}{2}y\right)\sinh\frac{y}{2} 
-\frac{1}{\sqrt{3}}\sin\left(\frac{\sqrt{3}}{2}y\right)\cosh\frac{y}{2}
\right]
\label{81g}
\end{eqnarray}
\begin{eqnarray}
\lefteqn{\sinh(ky) =\frac{1}{3}\sinh y
-\frac{2}{3}\cos\left(\frac{\sqrt{3}}{2}y\right)\sinh\frac{y}{2}
}\nonumber\\ 
&&+h\left[\frac{1}{3}\sinh y
+\frac{1}{3}\cos\left(\frac{\sqrt{3}}{2}y\right)\sinh\frac{y}{2} 
-\frac{1}{\sqrt{3}}\sin\left(\frac{\sqrt{3}}{2}y\right)\cosh\frac{y}{2}
\nonumber\right] \\
&&+k\left[\frac{1}{3}\sinh y
+\frac{1}{3}\cos\left(\frac{\sqrt{3}}{2}y\right)\sinh\frac{y}{2} 
+\frac{1}{\sqrt{3}}\sin\left(\frac{\sqrt{3}}{2}y\right)\cosh\frac{y}{2}
\right]
\label{81h}
\end{eqnarray}
The hyperbolic functions of a tricomplex number $x+h y+k z$ can then
be expressed in terms of the elementary functions with the aid of the addition
theorems Eqs. (\ref{81c}), (\ref{81d}) and of the expressions in  Eqs. 
(\ref{81e})-(\ref{81h}).

\section{Tricomplex power series}

A tricomplex series is an infinite sum of the form
\begin{equation}
a_0+a_1+a_2+\cdots+a_l+\cdots , 
\label{82}
\end{equation}
where the coefficients $a_n$ are tricomplex numbers. The convergence of 
the series (\ref{82}) can be defined in terms of the convergence of its 3 real
components. The convergence of a tricomplex series can however be studied
using tricomplex variables. The main criterion for absolute convergence 
remains the comparison theorem, but this requires a number of inequalities
which will be discussed further.

The modulus of a tricomplex number $u=x+hy+kz$ can be defined as
\begin{equation}
|u|=(x^2+y^2+z^2)^{1/2} .
\label{83}
\end{equation}
Since $|x|\leq |u|, |y|\leq |u|, |z|\leq |u|$, a property of absolute 
convergence established via a comparison theorem based on the modulus of the
series (\ref{82}) will ensure the absolute convergence of each real component
of that series.

The modulus of the sum $u_1+u_2$ of the tricomplex numbers $u_1, u_2$ fulfils
the inequality
\begin{equation}
||u_1|-|u_2||\leq |u_1+u_2|\leq |u_1|+|u_2| .
\label{84}
\end{equation}
For the product the relation is 
\begin{equation}
|u_1u_2|\leq \sqrt{3}|u_1||u_2| ,
\label{85}
\end{equation}
which replaces the relation of equality extant for regular complex numbers.
The equality in Eq. (\ref{85}) takes place for $x_1=y_1=z_1$ and 
$x_2=y_2=z_2$, i.e
when both tricomplex numbers lie on the trisector line $(t)$.
Using Eq. (\ref{65}), the relation (\ref{85}) can be written equivalently as 
\begin{equation}
\frac{2}{3}\delta_1^2\delta_2^2+\frac{1}{3}\sigma_1^2\sigma_2^2\leq 
3\left(\frac{2}{3}\delta_1^2+\frac{1}{3}\sigma_1^2\right)
\left(\frac{2}{3}\delta_2^2+\frac{1}{3}\sigma_2^2\right) ,
\label{85a}
\end{equation}
where $\delta_j^2=x_j^2+y_j^2+z_j^2-x_jy_j-x_jz_j-y_jz_j, \sigma_j=x_j+y_j+z_j,
j=1,2$, the equality taking place for $\delta_1=0,\delta_2=0$.
A particular form of Eq. (\ref{85}) is
\begin{equation}
|u^2|\leq \sqrt{3} |u|^2 ,
\label{86}
\end{equation}
and it can be shown that
\begin{equation}
|u^l|\leq 3^{(l-1)/2}|u|^l ,
\label{87}
\end{equation}
the equality in Eqs. (\ref{86}) and (\ref{87}) taking place for $x=y=z$.
It can be shown from Eq. (\ref{67}) that
\begin{equation}
|u^l|^2=\frac{2}{3}\delta^{2l}+\frac{1}{3}\sigma^{2l} ,
\label{87a}
\end{equation}
where $\delta^2=x^2+y^2+z^2-xy-xz-yz, \sigma=x+y+z$.
Then Eq. (\ref{87}) can also be written as
\begin{equation}
\frac{2}{3}\delta^{2l}+\frac{1}{3}\sigma^{2l}\leq
3^{l-1}\left(\frac{2}{3}\delta^2+\frac{1}{3}\sigma^2\right)^l ,
\label{87b}
\end{equation}
the equality taking place for $\delta=0$.
From Eqs. (\ref{85}) and (\ref{87}) it results that
\begin{equation}
|au^l|\leq 3^{l/2} |a| |u|^l .
\label{88}
\end{equation}
It can also be shown that
\begin{equation}
\left|\frac{1}{u}\right|\geq\frac{1}{|u|} ,
\label{89}
\end{equation}
the equality taking place for $\sigma^2=\delta^2$, or $xy+xz+yz=0$.

A power series of the tricomplex variable $u$ is a series of the form
\begin{equation}
a_0+a_1 u + a_2 u^2+\cdots +a_l u^l+\cdots .
\label{90}
\end{equation}
Since
\begin{equation}
\left|\sum_{l=0}^\infty a_l u^l\right| \leq \sum_{l=0}^\infty
3^{l/2}|a_l||u|^l ,
\label{91}
\end{equation}
a sufficient condition for the absolute convergence of this series is that
\begin{equation}
\lim_{n\rightarrow \infty}
\frac{\sqrt{3}|a_{l+1}||u|}{|a_l|}<1 .
\label{92}
\end{equation}
Thus the series is absolutely convergent for 
\begin{equation}
|u|<c_0,
\label{new86}
\end{equation}
where 
\begin{equation}
c_0=\lim_{l\rightarrow\infty} \frac{|a_l|}{\sqrt{3}|a_{l+1}|} .
\label{new87}
\end{equation}

The convergence of the series (\ref{90}) can be also studied with the aid of
a transformation which explicits the transverse and longitudinal parts
of the variable $u$ and of the constants $a_l$, 
\begin{equation}
x+hy+kz=v_1 e_1+\tilde v_1 \tilde e_1+v_+ e_+ , 
\label{116}
\end{equation}
where
\begin{equation}
v_1=\frac{2x-y-z}{2}, \:\:\tilde v_1=\frac{\sqrt{3}}{2}(y-z), \:\: v_+=x+y+z,
\label{117a}
\end{equation}
and
\begin{equation}
e_1=\frac{2-h-k}{3},\:\:\tilde e_1=\frac{h-k}{\sqrt{3}},
\:\:e_+=\frac{1+h+k}{3} .
\label{117b}
\end{equation}
The variables $v_1, \tilde v_1, v_+$ will be called the tricomplex
canonical variables,
and $e_1, \tilde e_1, e_+$ will be called the tricomplex canonical base.
In the geometric representation of Fig. 6, 
$e_1, \tilde e_1$ are situated in the plane $\Pi$, and $e_+$ is lying on the
trisector 
line $(t)$.  It can be checked that
\begin{equation}
e_1^2=e_1, \:\:\tilde e_1^2=-e_1,\:\: e_1\tilde e_1=\tilde e_1,\:\: e_1e_+=0,
\:\:\tilde e_1e_+=0,  
\:\:e_+^2=e_+. 
\label{118}
\end{equation}
The moduli of the bases in Eq. (\ref{118}) are
\begin{equation}
|e_1|=\sqrt{\frac{2}{3}},\;|\tilde e_1|
=\sqrt{\frac{2}{3}},\;|e_+|=\sqrt{\frac{1}{3}},
\label{118b}
\end{equation}
and it can be checked that
\begin{equation}
|x+hy+kz|^2=\frac{2}{3}(v_1^2+\tilde v_1^2)+\frac{1}{3}v_+^2 .
\label{118c}
\end{equation}

If $u=u^\prime u^{\prime\prime}$, the transverse and longitudinal
components are related by
\begin{equation}
v_1=v_1^\prime v_1^{\prime\prime}-\tilde v_1^\prime \tilde v_1^{\prime\prime},
\:\: 
\tilde v_1=v_1^\prime \tilde v_1^{\prime\prime}+\tilde v_1^\prime
v_1^{\prime\prime}, \:\: 
v_+=v_+^\prime v_+^{\prime\prime} ,
\label{119}
\end{equation}
which show that, upon multiplication, the transverse components obey the same
rules as the real and imaginary components of usual, two-dimensional complex
numbers, and the rule for the longitudinal component is that of the regular
multiplication of numbers.

If the constants in Eq. (\ref{90}) are $a_l=p_l+hq_l+kr_l$, and
\begin{equation}
a_{l1}=\frac{2p_l-q_l-r_l}{2}, \:\:\tilde a_{l1}=\frac{\sqrt{3}}{2}(q_l-r_l), 
\:\: a_{l+}=p_l+q_l+r_l,
\label{121}
\end{equation}
where $p_0=1, q_0=0, r_0=0$,
the series (\ref{90}) can be written as
\begin{equation}
\sum_{l=0}^\infty \left[a_{l1}e_1+\tilde a_{l1}\tilde e_1)(v_1e_1+\tilde
v_1\tilde e_1)^l 
+e_+ a_{l+} v_+^l\right].
\label{121b}
\end{equation}
The series in Eq. (\ref{121b}) is absolutely convergent for 
\begin{equation}
|v_+|<c_+,\:
(v_1^2+\tilde v_1^2)^{1/2}<c_1,
\label{n90}
\end{equation}
where 
\begin{equation}
c_+=\lim_{l\rightarrow\infty} \frac{|a_{l+}|}{|a_{l+1,u}|} ,\:
c_1=\lim_{l\rightarrow\infty} \frac
{\left(a_{l1}^2+ \tilde a_{l1}^2\right)^{1/2}}
{\left(a_{l+1,1}^2+ a_{l+1,2}^2\right)^{1/2}} .
\label{n91}
\end{equation}
The relations (\ref{n90}) and (\ref{118c}) show that the region of convergence
of the series (\ref{90}) is a cylinder
of radius $c_1\sqrt{2/3}$ and height $2c_+/\sqrt{3}$, having the
trisector line $(t)$ as axis and the origin as center, which can be called
cylinder of convergence, as shown in Fig. 7.

It can be shown that $c_1=(1/\sqrt{3})\;{\rm
min}(c_+,c_1)$, where ${\rm min}$ designates the smallest of
the numbers $c_+,c_1$. Using the expression of $|u|$ in
Eq. (\ref{118}), it can be seen that the spherical region of
convergence defined in Eqs. (\ref{new86}), (\ref{new87}) is a subset of the
cylindrical region of convergence defined in Eqs. (\ref{n90}) and (\ref{n91}).

\section{Analytic functions of tricomplex variables}

The derivative  
of a function $f(u)$ of the tricomplex variables $u$ is
defined as a function $f^\prime (u)$ having the property that
\begin{equation}
|f(u)-f(u_0)-f^\prime (u_0)(u-u_0)|\rightarrow 0 \:\:{\rm as} 
\:\:|u-u_0|\rightarrow 0 . 
\label{gs88}
\end{equation}
If the difference $u-u_0$ is not parallel to one of the nodal hypersurfaces,
the definition in Eq. (\ref{gs88}) can also 
be written as
\begin{equation}
f^\prime (u_0)=\lim_{u\rightarrow u_0}\frac{f(u)-f(u_0)}{u-u_0} .
\label{gs89}
\end{equation}
The derivative of the function $f(u)=u^m $, with $m$ an integer, 
is $f^\prime (u)=mu^{m-1}$, as can be seen by developing $u^m=[u_0+(u-u_0)]^m$
as
\begin{equation}
u^m=\sum_{p=0}^{m}\frac{m!}{p!(m-p)!}u_0^{m-p}(u-u_0)^p,
\label{gs90}
\end{equation}
and using the definition (\ref{gs88}).

If the function $f^\prime (u)$ defined in Eq. (\ref{gs88}) is independent of
the 
direction in space along which $u$ is approaching $u_0$, the function $f(u)$ 
is said to be analytic, analogously to the case of functions of regular complex
variables. \cite{3} 
The function $u^m$, with $m$ an integer, 
of the tricomplex variable $u$ is analytic, because the
difference $u^m-u_0^m$ is always proportional to $u-u_0$, as can be seen from
Eq. (\ref{gs90}). Then series of
integer powers of $u$ will also be analytic functions of the tricomplex
variable $u$, and this result holds in fact for any commutative algebra. 

If an analytic function is defined by a series around a certain point, for
example $u=0$, as
\begin{equation}
f(u)=\sum_{k=0}^\infty a_k u^k ,
\label{gs91a}
\end{equation}
an expansion of $f(u)$ around a different point $a$,
\begin{equation}
f(u)=\sum_{k=0}^\infty c_k (u-a)^k ,
\label{gs91aa}
\end{equation}
can be obtained by
substituting in Eq. (\ref{gs91a}) the expression of $u^k$ according to Eq.
(\ref{gs90}). Assuming that the series are absolutely convergent so that the
order of the terms can be modified and ordering the terms in the resulting
expression according to the increasing powers of $u-a$ yields
\begin{equation}
f(u)=\sum_{k,l=0}^\infty \frac{(k+l)!}{k!l!}a_{k+l} a^l (u-a)^k .
\label{gs91b}
\end{equation}
Since the derivative of order $k$ at $u=a$ of the function $f(u)$ , Eq.
(\ref{gs91a}), is 
\begin{equation}
f^{(k)}(a)=\sum_{l=0}^\infty \frac{(k+l)!}{l!}a_{k+l} a^l ,
\label{gs91c}
\end{equation}
the expansion of $f(u)$ around $u=a$, Eq. (\ref{gs91b}), becomes
\begin{equation}
f(u)=\sum_{k=0}^\infty \frac{1}{k!} f^{(k)}(a)(u-a)^k ,
\label{gs91d}
\end{equation}
which has the same form as the series expansion of the usual
2-dimensional complex functions. 
The relation (\ref{gs91d}) shows that the coefficients in the series expansion,
Eq. (\ref{gs91aa}), are
\begin{equation}
c_k=\frac{1}{k!}f^{(k)}(a) .
\label{gs92}
\end{equation}

The rules for obtaining the derivatives and the integrals of the basic
functions can 
be obtained from the series of definitions and, as long as these series
expansions have the same form as the corresponding series for the
2-dimensional complex functions, the rules of derivation and integration remain
unchanged. 

If the tricomplex function $f(u)$
of the tricomplex variable $u$ is written in terms of 
the real functions $F(x,y,z),G(x,y,z),H(x,y,z)$ of real variables $x,y,z$ as
\begin{equation}
f(u)=F(x,y,z)+hG(x,y,z)+kH(x,y,z), 
\label{98d}
\end{equation}
then relations of equality 
exist between partial derivatives of the functions $F,G,H$. These relations can
be obtained by writing the derivative of the function $f$ as
\begin{eqnarray}
\lefteqn{\frac{1}{\Delta x+h\Delta y +k\Delta z} 
\left[\frac{\partial F}{\partial x}\Delta x+
\frac{\partial F}{\partial y}\Delta y+
\frac{\partial F}{\partial z}\Delta z
+h\left(\frac{\partial G}{\partial x}\Delta x+
\frac{\partial G}{\partial y}\Delta y+
\frac{\partial G}{\partial z}\Delta z\right) \right.}\nonumber\\
&&\left. 
+k\left(\frac{\partial H}{\partial x}\Delta x+
\frac{\partial H}{\partial y}\Delta y+
\frac{\partial H}{\partial z}\Delta z\right)\right] ,
\label{99}
\end{eqnarray}
where the difference appearing in Eq. {(98)} is
$u-u_0=\Delta x+h\Delta y +k\Delta z$. The relations between 
the partials 
derivatives of the functions $F, G, H$ are obtained by setting successively in 
Eq. (\ref{99}) $\Delta x\rightarrow 0, \Delta y=0,  \Delta z=0$; then
$\Delta x= 0, \Delta y\rightarrow 0,  \Delta z=0$; and 
$\Delta x=0, \Delta y=0,  \Delta z\rightarrow 0$. The relations are
\begin{equation}
\frac{\partial F}{\partial x} = \frac{\partial G}{\partial y},\:\:
\frac{\partial G}{\partial x} = \frac{\partial H}{\partial y},\:\:
\frac{\partial H}{\partial x} = \frac{\partial F}{\partial y},
\label{100}
\end{equation}
\begin{equation}
\frac{\partial F}{\partial x} = \frac{\partial H}{\partial z},\:\:
\frac{\partial G}{\partial x} = \frac{\partial F}{\partial z},\:\:
\frac{\partial H}{\partial x} = \frac{\partial G}{\partial z},
\label{101}
\end{equation}
\begin{equation}
\frac{\partial G}{\partial y} = \frac{\partial H}{\partial z},\:\:
\frac{\partial H}{\partial y} = \frac{\partial F}{\partial z},\:\:
\frac{\partial F}{\partial y} = \frac{\partial G}{\partial z} .
\label{102}
\end{equation}
The relations (\ref{100})-(\ref{102}) are analogous to the Riemann relations
for the real and imaginary components of a complex function. It can be shown
from Eqs. (\ref{100})-(\ref{102}) that the components $F$ solutions of the
equations 
\begin{equation}
\frac{\partial^2 F}{\partial x^2}-\frac{\partial^2 F}{\partial y\partial z}=0,
\:\: 
\frac{\partial^2 F}{\partial y^2}-\frac{\partial^2 F}{\partial x\partial z}=0,
\:\: 
\frac{\partial^2 F}{\partial z^2}-\frac{\partial^2 F}{\partial x\partial y}=0,
\label{103}
\end{equation}
\begin{equation}
\frac{\partial^2 G}{\partial x^2}-\frac{\partial^2 G}{\partial y\partial z}=0,
\:\: 
\frac{\partial^2 G}{\partial y^2}-\frac{\partial^2 G}{\partial x\partial z}=0,
\:\: 
\frac{\partial^2 G}{\partial z^2}-\frac{\partial^2 G}{\partial x\partial y}=0,
\label{104}
\end{equation}
\begin{equation}
\frac{\partial^2 H}{\partial x^2}-\frac{\partial^2 H}{\partial y\partial z}=0,
\:\: 
\frac{\partial^2 H}{\partial y^2}-\frac{\partial^2 H}{\partial x\partial z}=0,
\:\: 
\frac{\partial^2 H}{\partial z^2}-\frac{\partial^2 H}{\partial x\partial y}=0.
\label{105}
\end{equation}
It can also be shown that the
differences $F-G, F-H, G-H$ are solutions of the equation of Laplace,
\begin{equation}
\Delta (F-G)=0,\:\:\Delta (F-H)=0,\:\:\Delta
(G-H)=0,\:\:\Delta=\frac{\partial^2}{\partial x^2}+\frac{\partial^2}{\partial
y^2}+\frac{\partial^2}{\partial z^2} 
\label{106}
\end{equation}

If a geometric transformation is considered in which to a point $u$ is
associated the point $f(u)$, 
it can be shown that the tricomplex function $f(u)$ transforms a straight line
parallel to the trisector line $(t)$ in a straight line parallel to $(t)$, and
transforms a plane parallel to the nodal plane $\Pi$ in a plane parallel to
$\Pi$. 
A transformation generated by a tricomplex function $f(u)$ does not conserve in
general the angle of intersecting lines.

\section{Integrals of tricomplex functions}

The singularities of tricomplex functions arise from terms of the form
$1/(u-a)^m$, with $m>0$. Functions containing such terms are singular not
only at $u=a$, but also at all points of a plane $(\Pi_a)$
through the point $a$ and
parallel to the nodal plane $\Pi$ and at all points of a
straight line $(t_a)$ 
passing through $a$ and parallel to the trisector line $(t)$.

The integral of a tricomplex function between two points $A, B$ along a path
situated in a region free of singularities is independent of path, which means
that the integral of an analytic function along a loop situated in a region
free from singularities is zero,
\begin{equation}
\oint_\Gamma f(u) du = 0,
\label{107}
\end{equation}
where it is supposed that a surface $S$ spanning the 
closed loop $\Gamma$ is not intersected by any of
the planes and is not threaded by any of the lines associated with the
singularities of the function $f(u)$. Using the expression, Eq. (\ref{98d})
for $f(u)$ and the fact that $du=dx+hdy+kdz$, the explicit form of the 
integral in Eq. (\ref{107}) is
\begin{equation}
\oint _\Gamma f(u) du = \oint_\Gamma
[Fdx+Hdy+Gdz+h(Gdx+Fdy+Hdz)+k(Hdx+Gdy+Fdz)] .
\label{108}
\end{equation}
If the functions $F, G, H$ are regular on a surface $S$
spanning the loop $\Gamma$,
the integral along the loop $\Gamma$ can be transformed with the aid of the
theorem of Stokes in an integral over the surface $S$ of terms of the form
$\partial H/\partial x -  \partial F/\partial y, \:\:
\partial G/\partial x - \partial F/\partial z,\:\:
\partial G/\partial y - \partial H/\partial z, \dots$
which are equal to zero by Eqs. (\ref{100})-(\ref{102}), and this proves Eq.
(\ref{107}). 

If there are singularities on the surface $S$, the integral $\oint f(u) du$
is not necessarily equal to zero. If $f(u)=1/(u-a)$ and the loop $\Gamma_a$ is
situated in the half-space above the plane $(\Pi_a)$ and encircles once
the line $(t_a)$, then 
\begin{equation}
\oint_{\Gamma_a} \frac{du}{u-a}=\frac{2\pi}{\sqrt{3}}(h-k) .
\label{109}
\end{equation}
This is due to the fact that the integral of $1/(u-a)$ 
along the loop $\Gamma_a$ is equal to the integral of $1/(u-a)$ along a circle
$(C_a)$ with the center on the line $(t_a)$ and perpendicular to this line,
as shown in Fig. 8. 
\begin{equation}
\oint_{\Gamma_a} \frac{du}{u-a}=\oint_{C_a} \frac{du}{u-a} ,
\label{110}
\end{equation}
this being a corrolary of Eq. (\ref{107}). The integral on the right-hand side
of Eq. (\ref{110}) can be evaluated with the aid of the trigonometric form 
Eq. (\ref{63}) of the tricomplex quantity $u-a$, so that
\begin{equation}
\frac{du}{u-a}=\frac{h-k}{\sqrt{3}}d\phi ,
\label{111}
\end{equation}
which by integration over $d\phi$ from 0 to $2\pi$ yields Eq. (\ref{109}).

The integral  $\oint_{\Gamma_a} du(u-a)^m$, 
with $m$ an integer number not equal to -1, is 
equal to zero, because $\int du (u-a)^m=(u-a)^{m+1}/(m+1)$, and 
$(u-a)^{m+1}/(m+1)$ is singlevalued,
\begin{equation}
\oint_{\Gamma_a}du(u-a)^m=0, \:\:{\rm for}\: m\: {\rm integer},
\: m\not=-1. 
\label{112}
\end{equation}
If $f(u)$ is an analytic tricomplex function which can be expanded in a
series as written in Eq. (\ref{gs91aa}), and the expansion holds on the curve
$\Gamma$ and on a surface spanning $\Gamma$, then from Eqs. (\ref{111}) and
(\ref{112}) it follows that
\begin{equation}
\oint_\Gamma \frac{f(u)du}{u-a}=\frac{2\pi}{\sqrt{3}}(h-k)f(a) .
\label{113}
\end{equation}
Substituting in the right-hand side of 
Eq. (\ref{113}) the expression of $f(u)$ in terms of the real 
components $F, G, H$,  Eq. (\ref{98d}), at $u=a$, yields
\begin{equation}
\oint_\Gamma \frac{f(u)du}{u-a}=\frac{2\pi}{\sqrt{3}}[H-G+h(F-H)+k(G-F)] .
\label{114}
\end{equation}
Since the sum of the real components in the paranthesis from the right-hand
side of Eq. (\ref{114}) is equal to zero, this equation
determines only the differences between the components $F, G, H$. 
If $f(u)$ can be expanded as written in Eq. (\ref{gs91aa}) on 
$\Gamma$ and on a surface spanning $\Gamma$, then from Eqs. (\ref{109}) and
(\ref{112}) it also results that
\begin{equation}
\oint_\Gamma \frac{f(u)du}{(u-a)^{m+1}}
=\frac{2\pi}{\sqrt{3}m!}(h-k)f^{(m)}(a) ,
\label{115}
\end{equation}
where the fact that has been used that the derivative $f^{(m)}(a)$ of order $m$
of $f(u)$ at $u=a$ is related to the expansion coefficient in Eq.
(\ref{gs91aa}) 
according to Eq. (\ref{gs92}). The relation (\ref{115}) can also be obtained by
successive derivations of Eq. (\ref{113}).

If a function $f(u)$ is expanded in positive and negative powers of $u-u_j$,
where $u_j$ are fourcomplex constants, $j$ being an index, the integral of $f$
on a closed loop $\Gamma$ is determined by the terms in the expansion of $f$
which are of the form $a_j/(u-u_j)$,
\begin{equation}
f(u)=\cdots+\sum_j\frac{a_j}{u-u_j}+\cdots .
\label{115b}
\end{equation}
In Eq. (\ref{115b}), $u_j$ is the pole and $a_j$ is the residue relative to the
pole $u_j$.
Then the integral of $f$ on a closed loop $\Gamma$ is
\begin{equation}
\oint_\Gamma f(u) du = 
\frac{2\pi}{\sqrt{3}}(h-k) \sum_j \:\;
{\rm int}(u_{j\Pi},\Gamma_{\Pi}) a_j ,
\label{115c}
\end{equation}
where the functional int($M,C$), defined for a point $M$ and a
closed curve $C$ in a two-dimensional plane, is given by 
\begin{equation}
{\rm int}(M,C)=\left\{
\begin{array}{l}
1 \;\:{\rm if} \;\:M \;\:{\rm is \;\:an \;\:interior \;\:point \;\:of} \;\:C
,\\  
0 \;\:{\rm if} \;\:M \;\:{\rm is \;\:exterior \;\:to}\:\; C\\
\end{array}\right. 
\label{115d}
\end{equation}
and $u_{j\Pi},\Gamma_{\Pi}$ are the projections of the point $u_j$ and
of the curve $\Gamma$ on the nodal plane $\Pi$, as shown in Fig. 9.

\section{Factorization of tricomplex polynomials}

A polynomial of degree $m$ of the tricomplex variable $u=x+hy+kz$ has the form
\begin{equation}
P_m(u)=u^m+a_1 u^{m-1}+\cdots+a_{m-1} u +a_m ,
\label{115e}
\end{equation}
where the constants are in general tricomplex numbers, $a_l=p_l+hq_l+kr_l$,
$l=1,\cdots,m$.  
In order to write the polynomial $P_m(u)$ as a product of factors, the variable
$u$ and the constants $a_l$ will be written in the form which explicits the
transverse and longitudinal components, 
\begin{equation}
P_m(u)=\sum_{l=0}^m (a_{l1}e_1+\tilde a_{l1}\tilde e_1)(v_1e_1
+\tilde v_1\tilde e_1)^{m-l}
+e_+\sum_{l=0}^m a_{l+} v_+^{m-l} ,
\label{120}
\end{equation}
where the constants have been defined previously in Eq. (\ref{121}).
Due to the properties in Eq. (\ref{118}), the transverse part of the
polynomial $P_m(u)$ can be written as a product of linear factors of the form
\begin{equation}
\sum_{l=0}^m (a_{l1}e_1+\tilde a_{l1}\tilde e_1)(v_1e_1+\tilde v_1\tilde e_1)^{m-l}=
\prod_{l=1}^m [(v_1-v_{l1})e_1+(\tilde v_1-\tilde v_{l1})\tilde e_1] ,
\label{122}
\end{equation}
where the quantities $v_{l1}, \tilde v_{l1}$ are real numbers.
The longitudinal part of $P_m(u)$, Eq. (\ref{120}), can be written as a product
of linear or quadratic factors with real coefficients, or as a product of
linear factors which, if imaginary, appear always in complex conjugate pairs.
Using the latter form for the simplicity of notations, the longitudinal part
can be written as
\begin{equation}
\sum_{l=0}^m a_{l+} v_+^{m-l}=\prod_{l=1}^m (v_+-v_{l+}) ,
\label{123}
\end{equation}
where the quantities $v_{l+}$ appear always in complex conjugate pairs.
Due to the orthogonality of the transverse and longitudinal components, Eq.
(\ref{118}), the polynomial $P_m(u)$ can be written as a product of factors of
the form  
\begin{equation}
P_m(u)=\prod_{l=1}^m [(v_1-v_{l1})e_1+(\tilde v_1-\tilde v_{l1})\tilde e_1
+ (v_+-v_{l+})e_+] .
\label{124}
\end{equation}
These relations can be written with the aid of Eqs. (\ref{116}) as
\begin{eqnarray}
P_m(u)=\prod_{l=1}^m (u-u_l),
\label{125}
\end{eqnarray}
where
\begin{eqnarray}
u_l=v_{l1}e_1+\tilde v_{l1}\tilde e_1+ v_{l+}e_+ .
\label{126}
\end{eqnarray}
The roots $v_{l+}$ and the roots $v_{l1}e_1+\tilde v_{l1}\tilde e_1$
defined in Eq. (\ref{122}) may be ordered arbitrarily.
This means that Eq. (\ref{126}) gives sets of $m$ roots
$u_1,...,u_m$ of the polynomial $P_m(u)$, 
corresponding to the various ways in which the roots 
$v_{l+}, v_{l1}e_1+\tilde v_{l1}\tilde e_1$ are ordered according to $l$ in
each group. Thus, while the tricomplex components in Eq. (\ref{120}) taken
separately have unique factorizations, the polynomial $P_m(q)$ can be written
in many different ways as a product of linear factors. 

If $P(u)=u^2-1$, the degree is $m=2$, the coefficients of the polynomial are
$a_1=0, a_2=-1$, the coefficients defined in Eq. (\ref{121}) are 
$a_{21}=-1, a_{22}=0, a_{2u}=-1$. 
The expression of $P(u)$, Eq. (\ref{120}), is  
$(e_1v_1+\tilde e_1\tilde v_1)^2-e_1+e_+ (v_+^2-1)$. The factorizations in Eqs.
(\ref{122}) 
and (\ref{123}) are
$(e_1v_1+\tilde e_1\tilde v_1)^2-e_1=[e_1(v_1+1)+\tilde e_1\tilde
v_1][e_1(v_1-1)+\tilde e_1\tilde v_1]$ and  
$v_+^2-1=(v_++1)(v_+-1)$. The factorization of $P(u)$, Eq. (\ref{125}), is
$P(u)=(u-u_1)(u-u_2)$, where according to Eq. (\ref{126}) the roots are
$u_1=\pm e_1\pm  e_+, u_2=-v_1$. If $e_1$ and $e_+$ are expressed with the
aid of Eq. (\ref{117b}) in terms of $h$ and $k$, the factorizations of $P(u)$
are obtained as
\begin{equation}
u^2-1=(u+1)(u-1) ,
\label{127}
\end{equation}
or as
\begin{equation}
u^2-1=\left(u+\frac{1-2h-2k}{3}\right) 
\left(u-\frac{1-2h-2k}{3}\right) .
\label{128}
\end{equation}
It can be checked that 
$(\pm e_1\pm  e_+)^2=e_1+e_+=1$.

\section{Representation of tricomplex complex numbers by irreducible
matrices} 

If the matrix in Eq. (\ref{25}) representing the tricomplex number $u$ is
called $U$, and 
\begin{equation}
T=\left(
\begin{array}{ccc}
\sqrt{\frac{2}{3}}     &  -\frac{1}{\sqrt{6}} & -\frac{1}{\sqrt{6}} \\
0                      &  \frac{1}{\sqrt{2}}   & -\frac{1}{\sqrt{2}} \\
\frac{1}{\sqrt{3}}     & \frac{1}{\sqrt{3}}   & \frac{1}{\sqrt{3}}   \\
\end{array}\right),
\label{129}
\end{equation}
which is the matrix appearing in Eq. (\ref{13}), it can be checked that 
\begin{equation}
TUT^{-1}=\left(
\begin{array}{ccc}
x-\frac{y+z}{2}          & \frac{\sqrt{3}}{2}(y-z) & 0  \\
-\frac{\sqrt{3}}{2}(y-z) &x-\frac{y+z}{2}          & 0  \\
  0                      & 0                       &x+y+z          
\end{array}\right).
\label{130}
\end{equation}
The relations for the variables $ x-(y+z)/2, (\sqrt{3}/2)(y-z)$ and 
$x+y+z$ for the multiplication of tricomplex numbers have been written in
Eqs. (\ref{18}), (\ref{20}) and (\ref{21}). The matrices
$T U T^{-1}$  provide an irreducible representation
\cite{4} of the tricomplex numbers $u=x+hy+kz$, in terms of matrices with real
coefficients.

\section{Conclusions}

The operations of addition and multiplication of the tricomplex numbers
introduced in this 
work have a simple geometric interpretation based on the amplitude $\rho$,
polar angle $\theta$ and azimuthal angle $\phi$. An exponential form exists for
the tricomplex numbers, and a trigonometric form exists involving the variables
$\rho, \theta $ and 
$\phi$. The tricomplex functions defined by series of powers are analytic, and
the partial derivatives of the components of the tricomplex functions are
closely related. The integrals of tricomplex functions are independent of path
in regions where the functions are regular. The fact that the exponential form
of the tricomplex numbers depends on the cyclic variable $\phi$ leads to the
concept of pole and residue for integrals on closed paths. The polynomials of
tricomplex variables can be written as products of linear or quadratic factors.

\newpage

FIGURE CAPTIONS\\

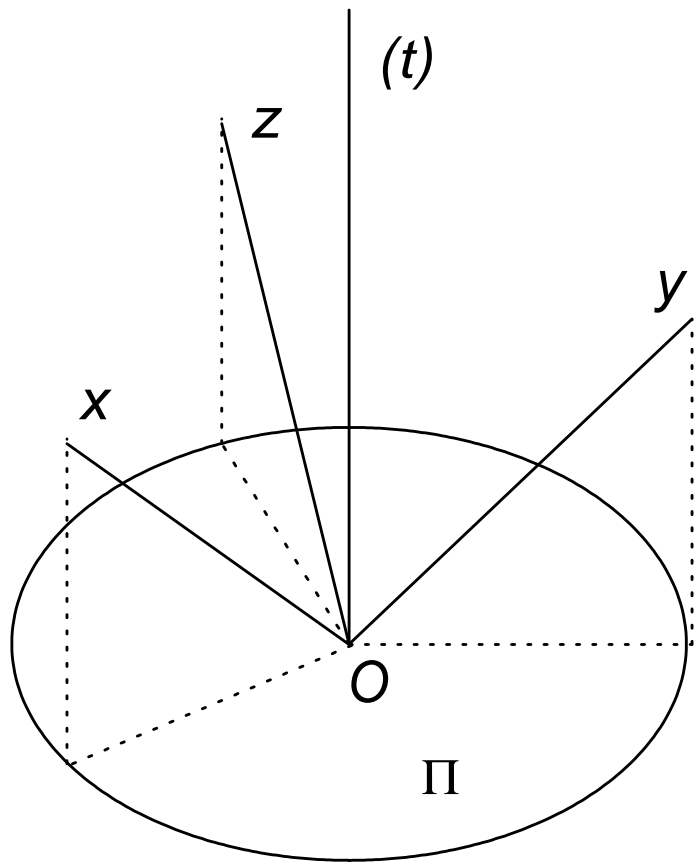
Fig. 1. Nodal plane $\Pi$, of equation $x+y+z=0$, and trisector line $(t)$, of
equation $x=y=z$, both passing through the origin $O$ of the rectangular axes
$x, y, z$.\\

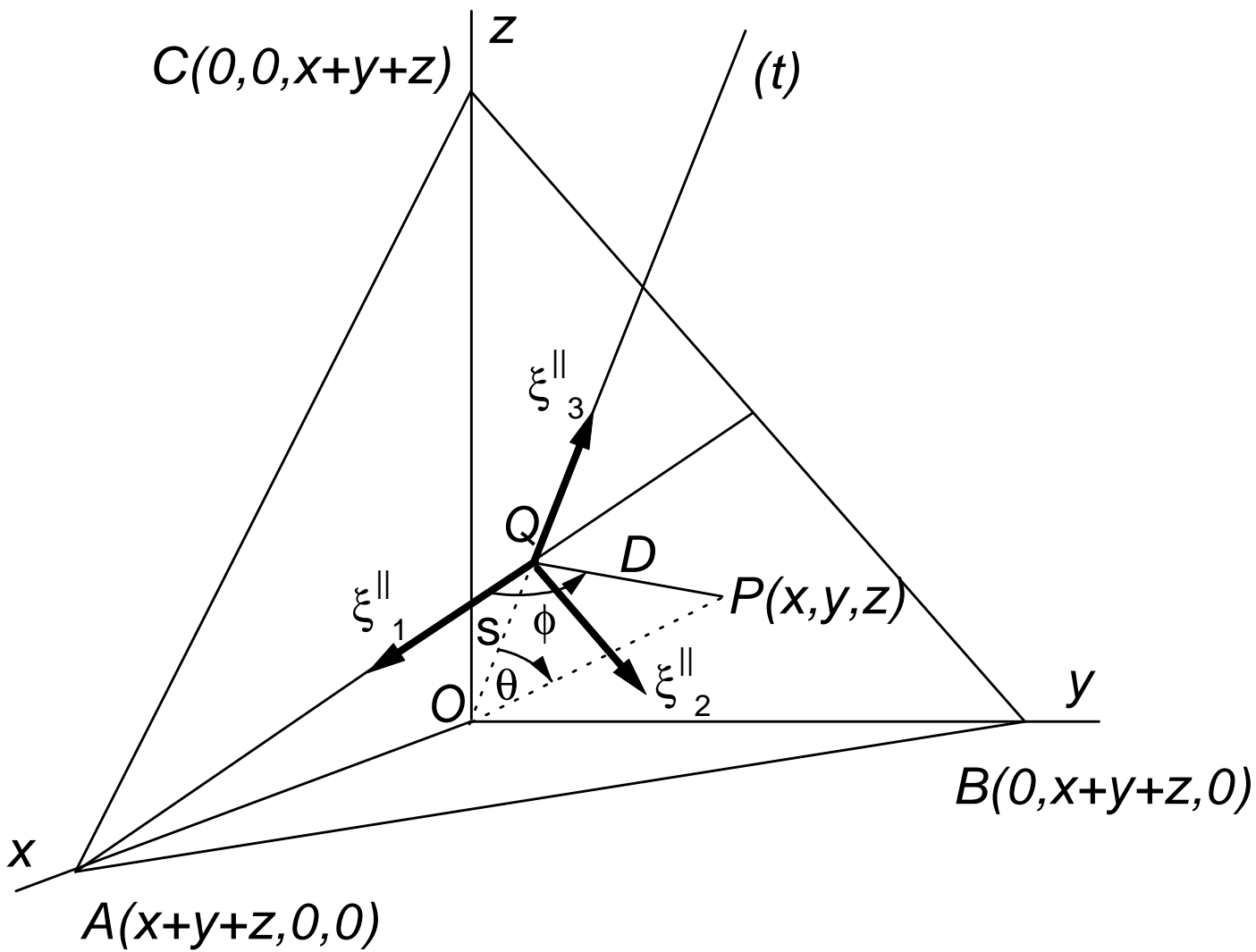
Fig. 2. Tricomplex variables $s, d, \theta, \phi$ for the tricomplex number
$x+hy+kz$, represented by the point $P(x,y,z)$. 
The azimuthal angle $\phi$ is shown in
in the plane parallel to $\Pi$, passing through $P$, which intersects 
the trisector line $(t)$ at $Q$ and the axes
of coordinates $x,y,z$ at the points $A, B, C$. The orthogonal axes
$\xi^\parallel _1,
\xi^\parallel_2, \xi^\parallel_3$ have the origin at $Q$.\\

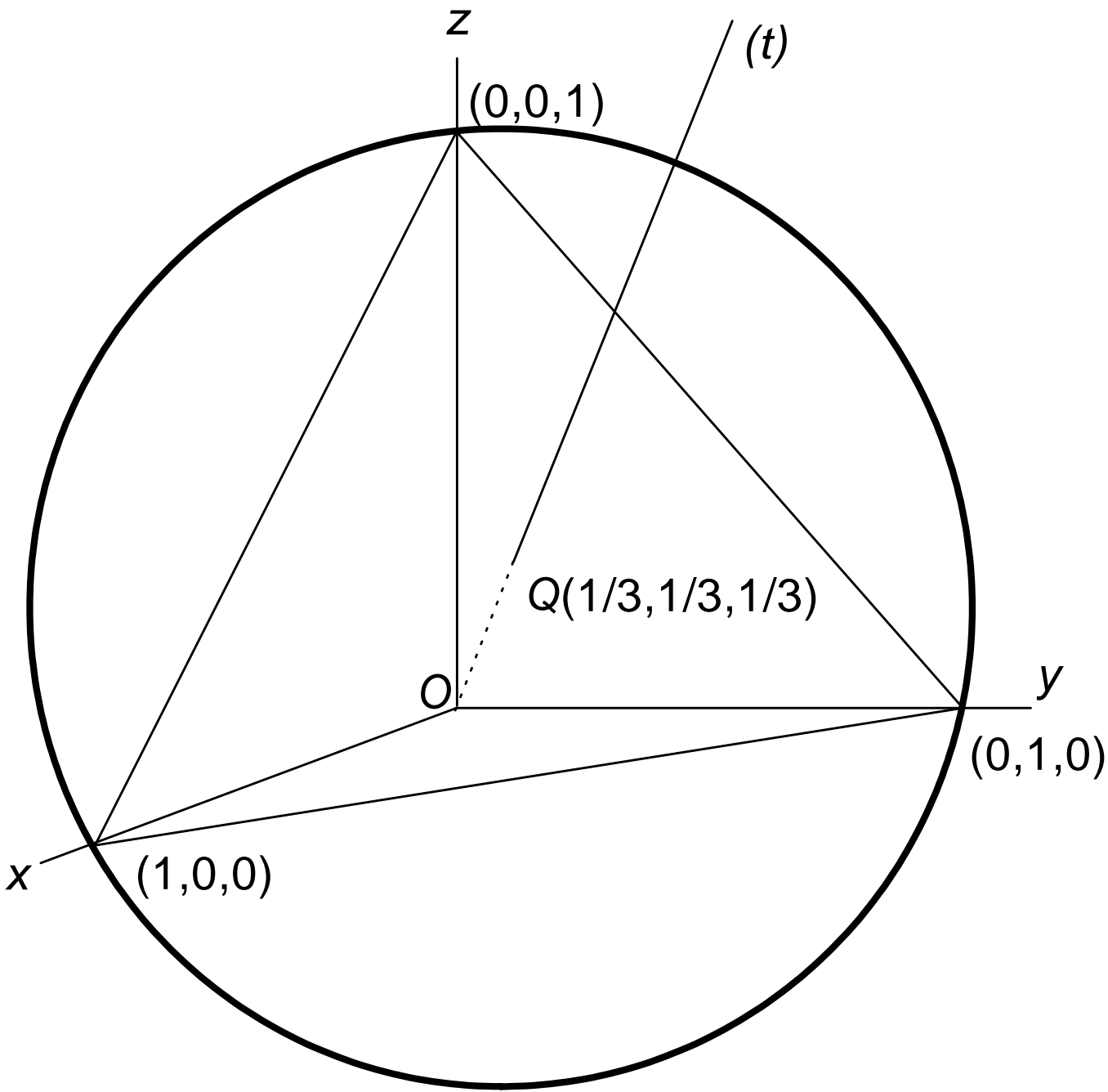
Fig. 3. Invariant circle for the multiplication of tricomplex numbers, lying
in a plane perpendicular to the trisector line and
passing through the points (1,0,0), (0,1,0) and (0,0,1). The center of the
circle is at the point $(1/3,1/3,1/3)$, and its radius is $\sqrt{2/3}$.\\

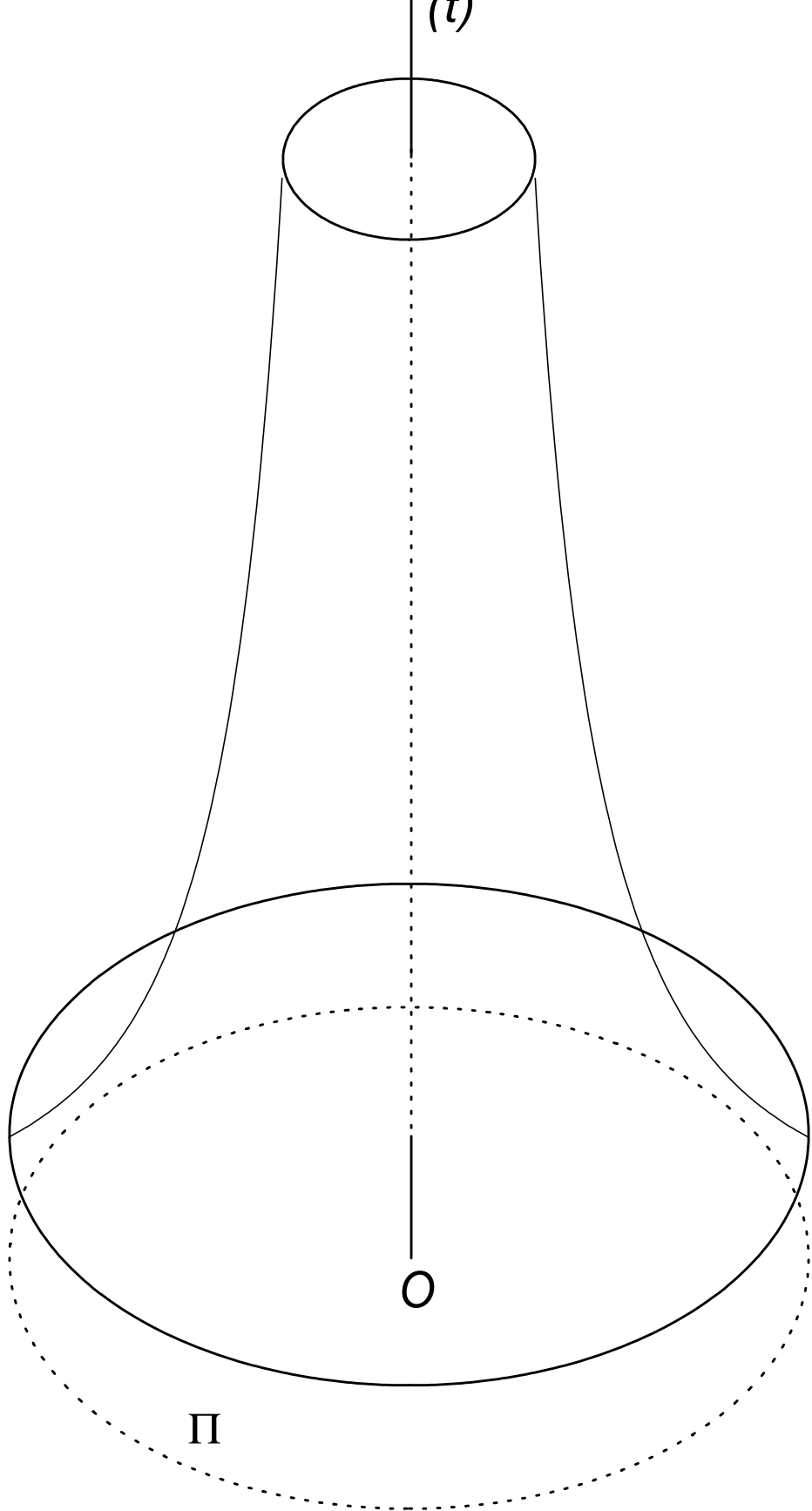
Fig. 4. Surfaces of constant $\rho$, which are surfaces of rotation
having the trisector line $(t)$ as axis.\\

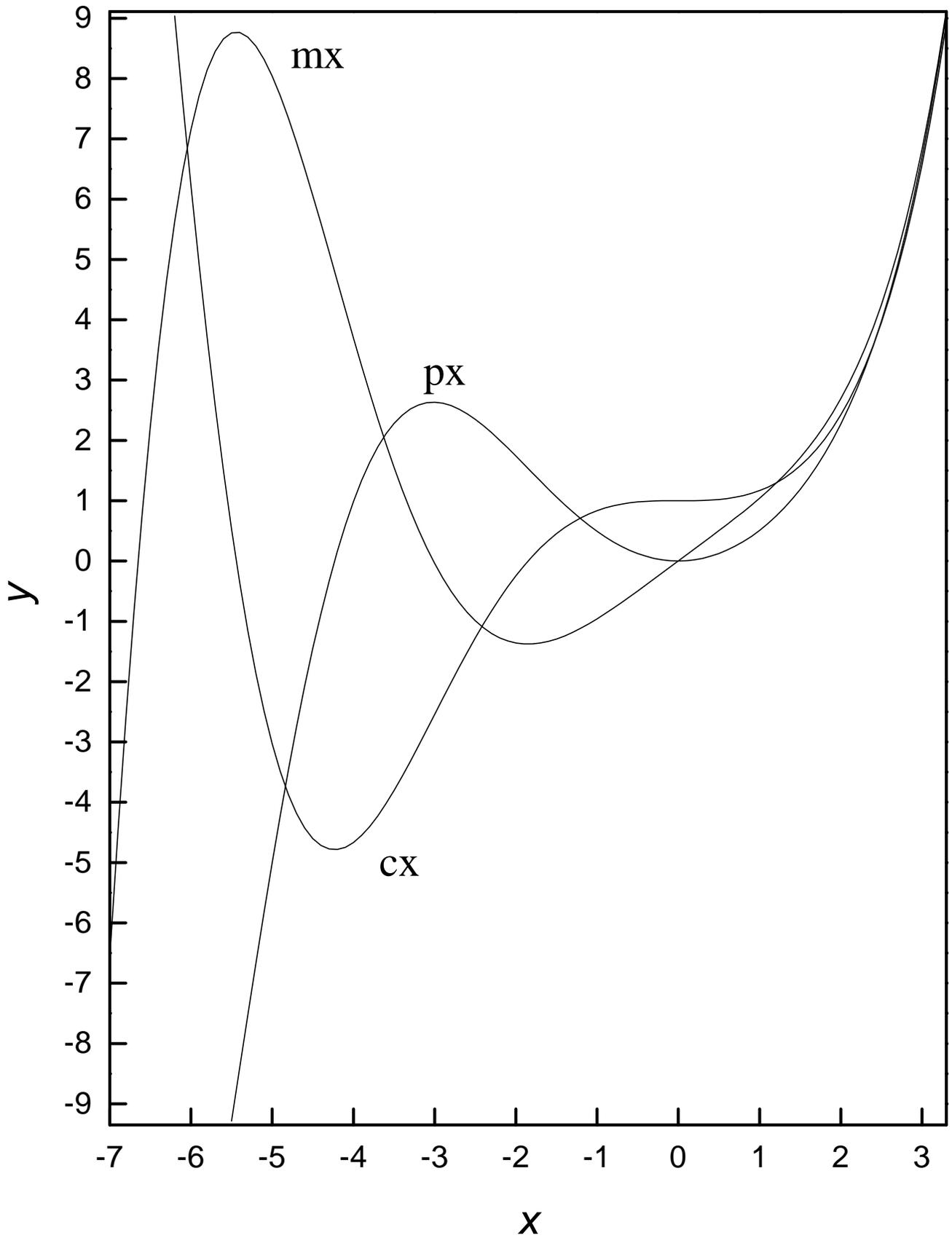
Fig. 5. Graphs of the cosexponential functions ${\rm cx}, {\rm mx},{\rm px}$.\\

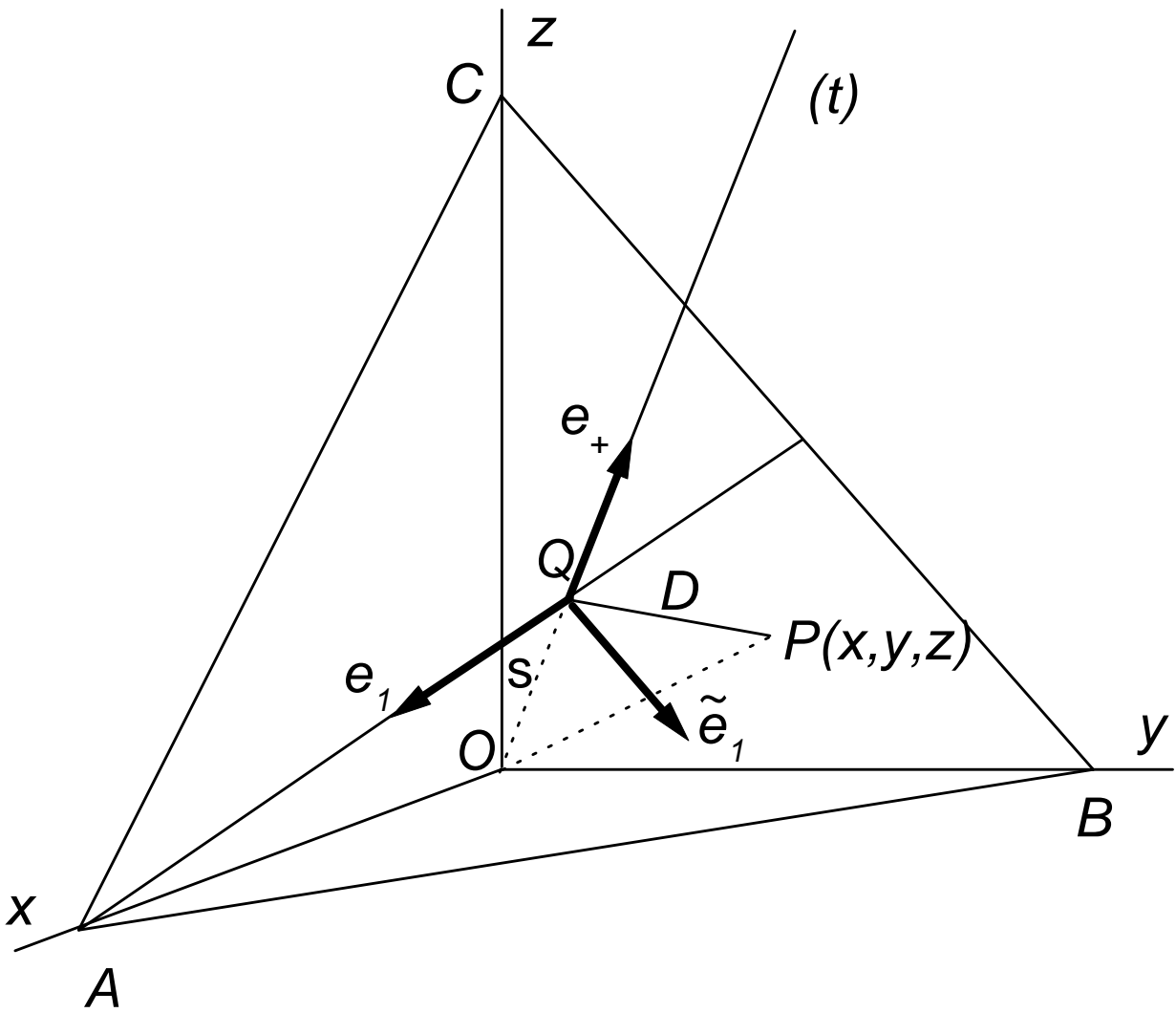
Fig. 6. Unit vectors $e_1, \tilde e_1, e_+$ of the orthogonal system of
coordinates with origin at $Q$.
The plane parallel to $\Pi$ passing through $P$ intersects 
the trisector line $(t)$ at $Q$ and the axes
of coordinates $x,y,z$ at the points $A, B, C$. \\

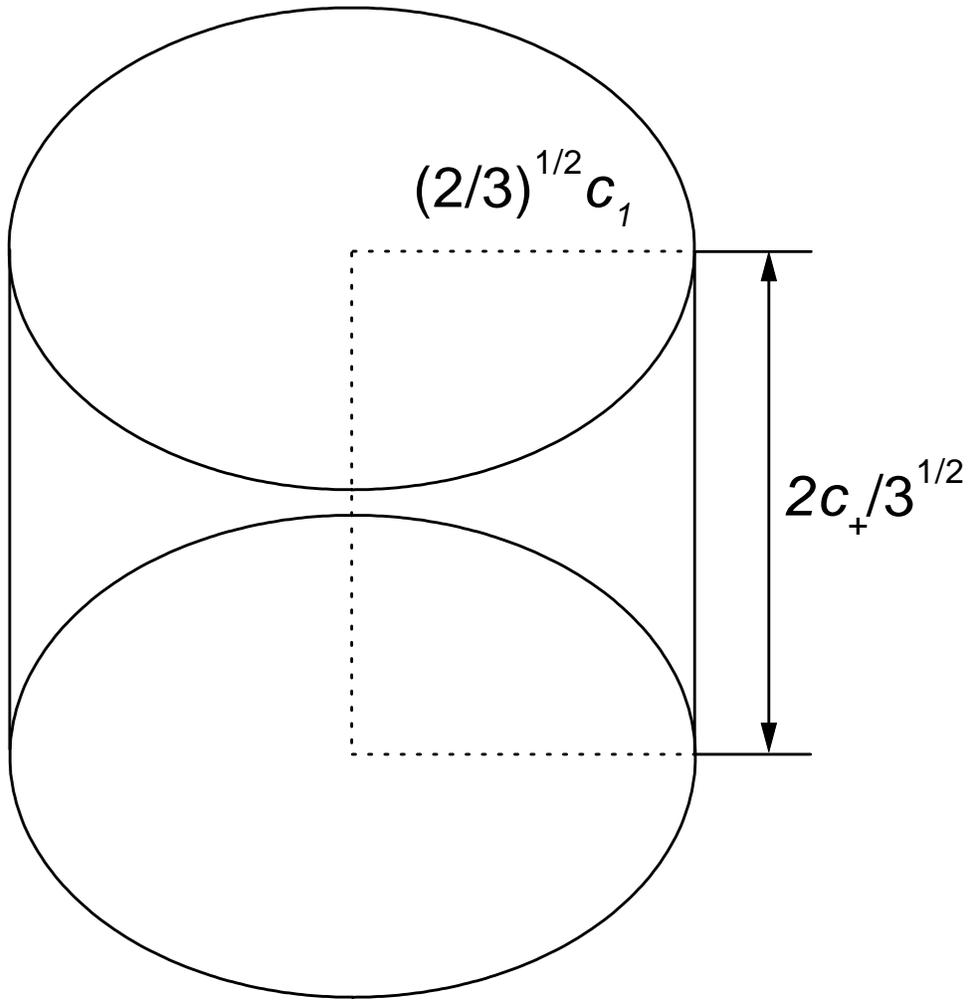
Fig. 7. Cylinder of convergence of tricomplex series, of radius $c_1\sqrt{2/3}$
and height $2c_+/\sqrt{3}$, having the axis parallel to the trisector line.\\

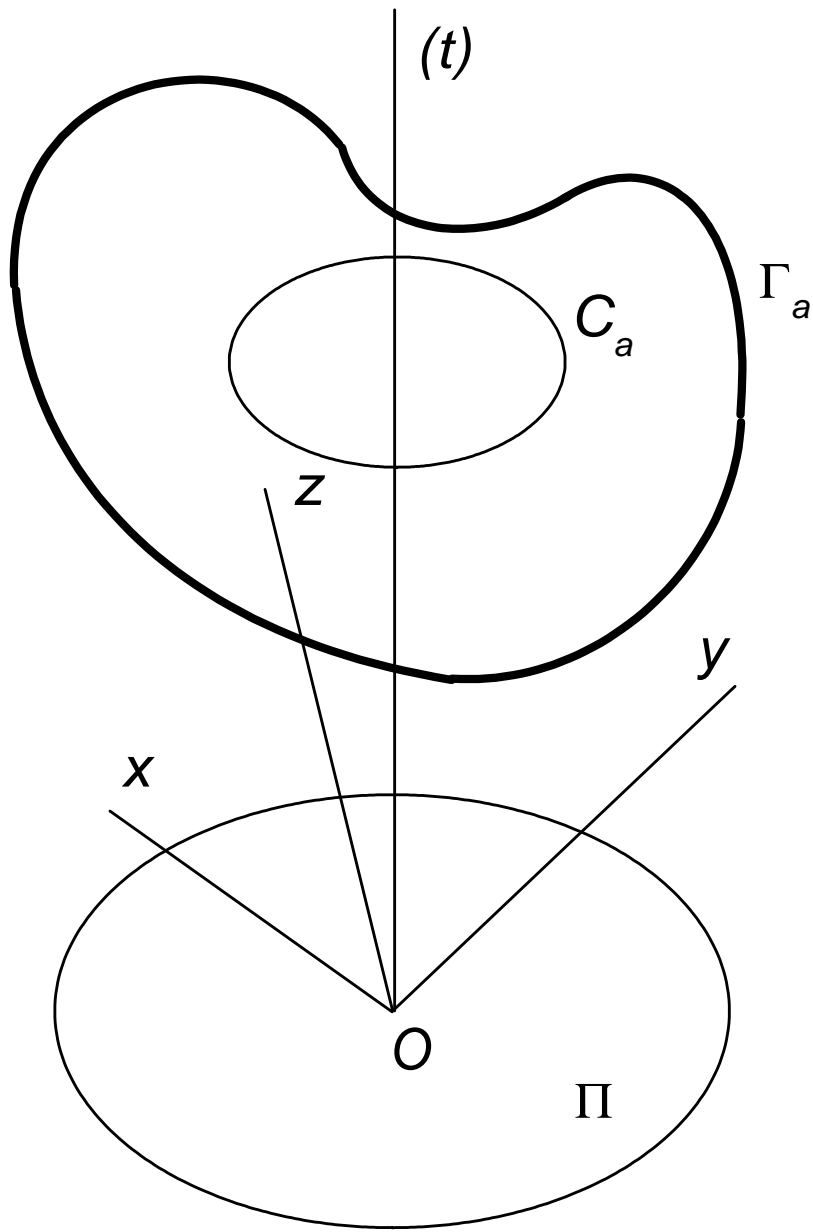
Fig. 8. The integral of $1/(u-a)$ 
along the loop $\Gamma_a$ is equal to the integral of $1/(u-a)$ along a circle
$(C_a)$ with the center on the line $(t_a)$ and perpendicular to this line.\\

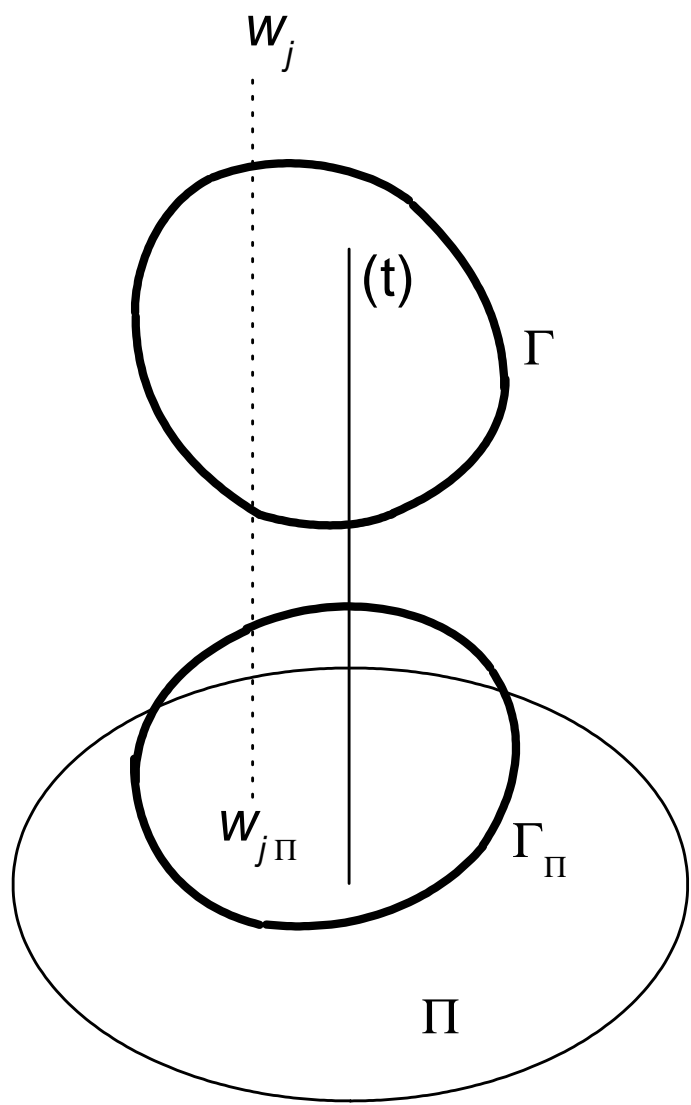
Fig. 9. Integration path $\Gamma$, pole $u_j$ and their projections
$\Gamma_{\Pi}, u_{j\Pi}$ on the nodal plane $\Pi$.\\

\end{document}